\documentclass[11pt,reqno]{amsart}
\usepackage{fullpage,amsfonts,amsmath,amscd,amssymb}
\input{xy}
\xyoption{all}
\theoremstyle{plain}

\newtheorem*{lem}{Lemma}
\newtheorem*{prop}{Proposition}
\newtheorem*{thm}{Theorem}
\newtheorem*{example}{Example}
\newtheorem*{cor}{Corollary}

\theoremstyle{remark}

\newtheorem*{rem}{Remark}

\newtheorem*{rems}{Remarks}

\newcommand{\x}[4]{\mbox{Ext}_{#1}^{#2}({#3},{#4})}
\newcommand{\hm}[3]{\mbox{Hom}_{#1}({#2},{#3})}
\newcommand{\cen}[1]{\mathfrak{c}_{\g}(#1)}
\newcommand{\im}{\text{Im}}
\newcommand{\ep}{\epsilon}
\newcommand{\g}{\mathfrak{g}}
\newcommand{\h}{\mathfrak{h}}
\newcommand{\n}{\mathfrak{n}}
\newcommand{\bp}{\mathfrak{b}^+}
\newcommand{\C}{\mathbb{C}}

\newcommand{\UE}{U_{\epsilon}(\mathfrak{g})}  
\newcommand{\UO}{U_0}  
\newcommand{\UGW}{U_{0}^{W\ltimes\Gamma}}
\newcommand{\UOG}{U_0^{\Gamma}}
\newcommand{\UL}{U_{0,\ell}}  
\newcommand{\ULG}{U_{0, \ell}^{\Gamma}}  
\newcommand{\ULGW}{U_{0, \ell}^{W\ltimes\Gamma}} 
\newcommand{\EO}{O_{\epsilon}[G]}   

\title[The ramification of centres]{The ramification of centres: Lie algebras in positive characteristic and quantised enveloping algebras.}
\author{Kenneth A. Brown}
\address{Mathematics Department, University of Glasgow, Glasgow G12 8QW}
\email{kab@maths.gla.ac.uk}
\author{Iain Gordon}
\address{Mathematics Department, J.C.M.B., King's Buildings, University of Edinburgh, Mayfield Road, EdinburghEH9 3JZ}
\email{igordon@maths.ed.ac.uk}
\thanks{ We are very grateful to J. Jantzen, L. Lebruyn, R. Marsh,
A. Premet and G. R\"{o}hrle for advice, discussions and
information. This research was begun while the second author was
studying for a Ph.D. at the University of Glasgow, with the support of
an EPSRC Research Studentship. The research of the first author was
partially supported by NATO Grant CRG 960250.
\newline
E-mail address for I.Gordon: igordon@maths.ed.ac.uk; fax number for I.Gordon: +44 131 650 6553.}
\begin{document}
\begin{abstract} 
Let $H$ be a Hopf algebra over the field $k$ which is a finite module
over a central affine sub-Hopf algebra $R$. Examples include
enveloping algebras $U(\g)$ of finite dimensional $k$-Lie algebras
$\g$ in positive characteristic and quantised enveloping algebras and
quantised function algebras at roots of unity. The ramification
behaviour of the maximal ideals of $Z(H)$ with respect to the
subalgebra $R$ is studied, and the conclusions are then applied to the
cases of classical and quantised enveloping algebras. In the case of
$U(\g)$ for $\g$ semisimple a conjecture of Humphreys \cite{hum3} on
the block structure of $U(\g)$ is confirmed. In the case of $U_{\ep}(\g)$ for
$\g$ semisimple and $\ep$ an odd root of unity we obtain a quantum analogue
of a result of Mirkovi\'{c} and Rumynin, \cite{mirrum}, and we fully describe the
factor algebras lying over the regular sheet, \cite{deckacpro007}. The blocks of $U_{\ep}(\g)$ are determined, and a necessary condition (which may also be sufficient) for a baby Verma $U_{\ep}(\g)$-module to be simple is obtained. 
\end{abstract}
\maketitle
\section{Introduction}
\subsection{}
\label{1.1}
Throughout $k$ will denote an algebraically closed field. In recent years common themes have become increasingly apparent in the representation theory of three important classes of $k$-algebras: the enveloping algebras $U(\g)$ of semisimple Lie algebras $\g$ in positive characteristic, the quantised enveloping algebras $\UE$ of semisimple Lie algebras at a root of unity $\epsilon$, and the quantised function algebras $\EO$ of semisimple groups $G$ at a root of unity $\epsilon$, \cite{jan6}, \cite{decpro3}, \cite{declyu1}. The common structure underlying these (and other related) classes is that of a triple
\begin{equation}
\label{triple}
R \subseteq Z \subseteq H
\end{equation}
of $k$-algebras, where $H$ is a Hopf algebra with centre $Z$, $Z$ being an affine domain, and $R$ is a sub-Hopf algebra of $H$, contained in $Z$, over which $H$ (and hence $Z$) are finite modules. The common strategy adopted in studying the (finite dimensional) representation theory of such an algebra is to study the finite dimensional $k$-algebras $H/\mathfrak{m}H$, as $\mathfrak{m}$ ranges across the maximal ideal spectrum of $R$.
\subsection{}
In this paper we continue the approach proposed and adopted in \cite{bro1}, \cite{brogoo} of looking for general results in the above setting which can then be interpreted and applied in the specific contexts mentioned above. Our starting point here is the following. Given a maximal ideal $\mathfrak{m}$ of $R$, how does the ramification behaviour of the maximal ideals of $Z$ lying over $\mathfrak{m}$ interact with the representation theory of $H/\mathfrak{m}H$? And how does this ramification behaviour vary as $\mathfrak{m}$ varies through $\mathrm{Maxspec}(R)$? We discuss these questions first in the abstract setting of a triple (\ref{triple}) in Section 2, and then consider classical and quantised enveloping algebras in Sections 3 and 4 respectively. (An analogous discussion for $\EO$, where more precise results can currently be proved than in the first two classes, is given in the sequel \cite{brogor2} to the present paper.)
\subsection{}
\label{1.3}
In Section 2, having first noted the easy fact that, \textit{in the setting} (1), \textit{the unramified locus of} Maxspec$(Z)$ \textit{is contained in the smooth locus}, we go on in Theorem \ref{azunram} to give a characterisation of an unramified point of Maxspec$(Z)$ under hypotheses which are satisfied in each of the three settings mentioned above. Thus, it is the main result of \cite{brogoo} that the smooth locus of Maxspec$(Z)$ coincides with the Azumaya locus of $H$ for each of the three classes listed in (\ref{1.1}); see Theorem \ref{az=sing}. (The \textit{Azumaya locus} of $H$ consists of those maximal ideals $M$ of $Z$ for which $H/MH$ is simple (artinian).) Theorem \ref{azunram} connects ramification with representation theory: it states that \textit{when the smooth locus of $Z$ coincides with the Azumaya locus a maximal ideal $M$ of $Z$ is unramified over $\mathfrak{m} = R \cap M$ if and only if $M$ is an Azumaya point and $H/MH$ is a projective $H/\mathfrak{m}H$-module.} Define a \textit{fully Azumaya point} $\mathfrak{m}$ of $R$ to be a maximal ideal $\mathfrak{m}$ of $R$ such that all the maximal ideals of $Z$ which lie over $\mathfrak{m}$ are in the Azumaya locus. Then we shall also be concerned to identify the \textit{fully Azumaya points} $\mathfrak{m}$ of $\mathrm{Maxspec}(R)$, and to describe the corresponding factors $H/\mathfrak{m}H$. \\
\indent The second theme of Section 2 is the problem of describing the \textit{blocks} of $H/\mathfrak{m}H$, for a maximal ideal $\mathfrak{m}$ of $R$. We point out in Proposition \ref{Blocks} that a ring-theoretic result of B.M\"{u}ller \cite{mul} applies here, to imply that
\begin{eqnarray*}
\textit{irreducible } H/\mathfrak{m}H-\textit{modules } V \textit{ and
} W \textit{ are in the same}\\ 
\textit{block of } H/\mathfrak{m}H \textit{ if and only if } \mathrm{Ann}_{Z}(V) = \mathrm{Ann}_{Z}(W).
\end{eqnarray*}
In particular, it follows at once that \textit{the number of blocks of $H/\mathfrak{m}H$ equals the number of maximal ideals of $Z$ which lie over $\mathfrak{m}$.}
\subsection{}
\label{lieintro}
In Section 3 we apply the abstract considerations of Section 2 to
the case where $H = U(\g)$ is the enveloping algebra of the Lie
algebra of a connected reductive group, $G$, whose derived subgroup is
simply-connected over a field $k$ of positive
characteristic $p$. Following \cite[6.3]{jan6} we assume throughout Section 3  
that (B) $p$ is odd and good for $\g$ and (C) the trace form on $\g$
is non-degenerate. See (\ref{az=sing}) for a discussion of these
hypotheses, whose key advantage is that if they are satisfied by $(\g,
p)$  then they also hold for $(\mathfrak{l}, p)$ for any Levi factor
$\mathfrak{l}$ of $\g$, so that one can rigourously employ the method of ``reduction to a nilpotent character" using the Morita theorem of Kac and Weisfeiler (Theorem \ref{morbabyV}). Thus one of our main aims in Section 3 is to formulate and prove in the above setting results previously known under stronger hypotheses on $\g$ and $p$. For example, \textit{we describe the centre $Z$ of $U(\g)$} (Theorem \ref{centreU}), weakening the hypotheses of earlier work of Veldkamp \cite{veld}, Mirkovi\'{c} and Rumynin\cite{mirrum}.\\
\indent Let $\g = \sum^{\oplus n}_{i=1}kx_{i}$. Then $R$ is the
$p$-\textit{centre} of $U(\g)$, a polynomial algebra in the variables
$\{x_i^p - x_i^{[p]} : 1 \leq i \leq n \}$, whose maximal ideals are
naturally parametrised by $\g^{*}$. Let $\mathfrak{m}_{\chi}$ be such
a maximal ideal and $U_{\chi} = U_{\chi}(\g):= U(\g)/\mathfrak{m}_{\chi}U(\g)$. The
starting point of our analysis is work of \cite{mirrum} giving a description (which we recall in \ref{mirrum}) of
$Z/\mathfrak{m}_{\chi}Z$ as a direct sum of local Frobenius
algebras. We use this description in Theorem \ref{Ugunram} and
Corollary \ref{babyproj} to characterise the unramified maximal ideals
of $Z$ lying over $\mathfrak{m}_{\chi}$ and to relate this to the
representation theory of $U_{\chi}$. For example we show there that
\textit{a maximal ideal of $Z$ lying over $\mathfrak{m}_{\chi}$ is
unramified if and only if the corresponding baby Verma module is a
projective $U_{\chi}$-module, and we give a condition equivalent to
this in the language of root system combinatorics of $\h^*$}. This
includes and extends earlier work of Friedlander and Parshall
\cite{frepar2}. We determine the fully Azumaya points of
$\mathrm{Maxspec}(R)$ - these are just the ideals
$\mathfrak{m}_{\chi}$ with $\chi$ regular (Proposition
\ref{azlocsec}); and we completely describe the structure of the
algebras $U_{\chi}$ in this case (Theorem \ref{regcase}), generalising
special cases obtained by Friedlander and Parshall, by Jantzen and by
Mirkovi\'{c} and Rumynin. We
also describe, for $\chi$ of standard Levi type, when the blocks of $U_{\chi}$ have finite
representation type, generalising a recent theorem of Nakano and
Pollack, \cite{nakpol}.  We emphasise again that our main purpose in this section has been to shed new light on old results by deriving them under uniform hypotheses and viewing them in terms of the relation between the central subalgebras $R$ and $Z$ of $U(\g)$. This applies also to our discussion of the blocks of $U_{\chi}$ in (\ref{lieblocks}).\\
\indent Let $T$ be a maximal torus of $G$, $W=N(T)/T$ the Weyl group of $G$
and $X=X(T)$, the character group of $T$. A basic problem in understanding
the representation theory of $U_{\chi}$ is to understand the
relationship between $Z/\mathfrak{m}_{\chi}Z$ and $Z(U_{\chi})$. For
generic (and in particular, for regular) $\chi$ these algebras are
isomorphic, but we record in (\ref{lieblocks}) a calculation due to
Premet which shows that \textit{the canonical homomorphism from
$Z/\mathfrak{m}_{\chi}Z$ to $Z(U_{\chi})$ is never onto when $\chi =
0$.} On the other hand, the result of M\"{u}ller mentioned in
(\ref{1.3}) implies that every idempotent of $Z(U_{\chi})$ is in the
image of $Z$. Hence we can determine the blocks of $U_{\chi}$, so
confirming a conjecture of J.E. Humphreys \cite[Section 18]{hum3} concerning the blocks of $U_{\chi}$ for $\g$ semisimple:
  \textit{For nilpotent $\chi$, the blocks of
$U_{\chi}$ are in bijection with the $W$-linkage classes in $X/pX$}.
\subsection{}
In Section 4 we study the case $H=U_{\ep}(\g)$, the quantised
enveloping algebra of the semisimple Lie algebra $\g=\text{Lie} (G)$ at a root of
unity, $\ep$. A significant problem in studying these algebras is the
behaviour of
the centraliser of a semisimple element of $G$. To control this we must
assume that $G$ is simply-connected and as a result that $U_{\ep}(\g)$
is the simply-connected form of the quantised enveloping
algebra. Under this restriction there is a description of the centre
of $U_{\ep}(\g)$ due to DeConcini, Kac and Procesi which is
qualitatively similar to the characteristic $p$ Lie algebra
situation. In particular $R$ is the ``$\ell$-centre'' of $U_{\ep}(\g)$
whose maximal ideals naturally correspond to elements of the ``big cell''
of $G$ (up to a finite covering). Given a maximal ideal
$\mathfrak{m}_{\chi}$ of $R$ let $U_{\chi}=
U_{\ep}(\g)/\mathfrak{m}_{\chi}U_{\ep}(\g)$. We prove (Theorem \ref{mainqgthm}) 
a quantum analogue of the result of Mirkovi\'{c} and Rumynin \cite{mirrum} describing
the primary components of $Z/\mathfrak{m_{\chi}}Z$. As in Section 3
we are able to use this to
\textit{characterise the unramified maximal ideals of $Z$ lying over
$\mathfrak{m}_{\chi}$} and give representation theoretic
consequences. We are also able to describe the blocks of
$U_{\chi}$ using M\"{u}ller's theorem. We give a \textit{necessary
condition for the simplicity of a quantised baby Verma module} in terms
of $\chi$ and root space combinatorics, allowing us to
completely describe the fully Azumaya points of $\text{Maxspec}(R)$
and also the structure of the corresponding algebra $U_{\chi}$. 
\section{Ramification}
\subsection{Index of ramification}
\label{indram}
Let $\theta: (R, \mathfrak{m}) \longrightarrow (Z, M)$ be a finite injective map 
of commutative local $k$-algebras having residue fields $k$. We define the 
\textit{index of ramification} of $R$ in $Z$ (along $\theta$) to be the least 
positive integer $i$ such that $M^i\subseteq \mathfrak{m}Z$. We say $R$ 
 is \textit{unramified} in $Z$ if the index of ramification equals one. 
Now let $\theta :R \longrightarrow Z$ be a finite injective map of affine 
commutative $k$-algebras and let $M$ be a maximal ideal of $Z$ and 
$\mathfrak{m}= M\cap R$. We define the index of ramification of $M$ (along 
$\theta$) to be the index of ramification of the finite embedding of local rings 
$\theta_M :R_{\mathfrak{m}} \longrightarrow Z_M$.
This defines a function
\[
i_{\theta}: \text{Maxspec}(Z) \longrightarrow \mathbb{N}.
\]
Let $\Omega_{Z/R}$ be the module of relative differentials (in other words the object 
which represents the functor $\text{Der}_R(Z, -)$) \cite[Chapter II.8]{har}).
\begin{lem}
Let $M$ be a maximal ideal of $Z$. Then $i_{\theta}(M) = 1$ if and only if 
$(\Omega_{Z/R})_M=0$.
\end{lem}
\begin{proof}
In view of \cite[Proposition 8.2A]{har} we can localise at $M$ and prove the 
analogous result for the local rings $R$ and $Z$. Let 
$\mathfrak{m} = M\cap R$ and let $\overline{Z} = Z/\mathfrak{m}Z$ with 
$\overline{M}$ the image 
of $M$ in $\overline{Z}$. Now, using \cite[p.186]{mat} and \cite[Chapter II, 
Proposition 8.7]{har}, we see 
\[
\frac{\Omega_{Z/R}}{M\Omega_{Z/R}} \cong \Omega_{Z/R} \otimes_Z \frac{Z}{M} 
 \cong \Omega_{\overline{Z}/k} \otimes_{\overline{Z}} k 
\cong  \frac{\overline{M}}{\overline{M}^2}.
\]
Since $\overline{Z}$ is local, Nakayama's lemma implies that $\overline{M} = 0$ 
if and only if $\Omega_{Z/R} =0$. 
\end{proof}

\subsection{The unramified locus}
\label{unloc}
For $\theta :R \longrightarrow Z$ as above we define the \textit{unramified 
locus (of $Z$ with respect to $R$)} to be the subset of unramified maximal 
ideals of $\mathrm{Maxspec}(Z)$. 
\begin{cor}
The unramified locus is an open set of $\mathrm{Maxspec}(Z)$, equal to the union 
of the unramified loci of the irreducible components of 
$\mathrm{Maxspec}(Z)$. 
\end{cor}
\begin{proof}
Let $\mathcal{M}$ be a coherent sheaf on a variety $X$ and suppose that 
$\mathcal{M}_x= 0$ for some $x\in X$. Then it follows immediately from 
Nakayama's lemma that there exists a neighbourhood $U$ of $x$ such that 
$\mathcal{M}|_U=0$. Since, by Lemma \ref{indram}, the unramified locus consists 
precisely of the points of $\mathrm{Maxspec}(Z)$ on which the (coherent) sheaf 
of 
relative differentials disappears, it is clear that the unramified locus is 
open.

Now let $M$ be a maximal ideal of $Z$ containing the minimal prime $P$. Then $M$ 
is unramified (along $\theta$) if and only if $M/P$ is unramified along the 
induced map
\[
\overline{\theta}: \frac{R}{P\cap R} \longrightarrow \frac{Z}{P}.
\]
Thus the unramified locus of $\text{Maxspec}(Z)$ is the union of the 
unramified loci of the irreducible components of $\text{Maxspec}(Z)$.
\end{proof}
\subsection{}
\label{nonemp}
In the cases of interest to us the unramified locus is non-empty, as the following lemma 
shows.
\begin{lem}
Suppose $\theta:R\longrightarrow Z$ is a finite injection of affine commutative 
integral domains. Then the unramified locus is 
non-empty if $\theta^*:\mathrm{Maxspec}(Z) \longrightarrow \mathrm{Maxspec}(R)$ 
is a 
separable morphism.
\end{lem}
\begin{proof}
Let $Q(R)$ and $Q(Z)$ be the quotient fields of $R$ and $Z$
respectively. Then, by definition, 
$\theta^*$ is separable if and only if $Q(Z)$ is a separable extension of 
$Q(R)$. If the latter hypothesis holds, then by \cite[Chapter II, Section 6, Theorem 4]{shaf} the
generic (set theoretic) fibre 
of $\theta^*$ contains $[Q(Z):Q(R)]$ elements. Since the $k$-dimension of 
$Z/\mathfrak{m}Z$ is $[Q(Z):Q(R)]$ for a generic maximal ideal $\mathfrak{m}$ of 
$R$, in the generic case the fibre is 
isomorphic to a direct product of copies of $k$ as required.
\end{proof}
\subsection{}
\label{smooth}
The following lemma will be crucial for us later when we want to compare 
representation theory with geometry.
\begin{lem}
Let $\theta :R \longrightarrow Z$ be a finite injection of affine commutative 
$k$-algebras and suppose that $\mathrm{Maxspec}(R)$ is smooth. Then the unramified 
locus of $\mathrm{Maxspec}(Z)$ is contained in the smooth locus.
\end{lem}
\begin{proof}
Suppose $M\in \text{Maxspec}(Z)$ is unramified and let $\mathfrak{m}=M\cap R$. 
Then
\[
\mathrm{Krull dim}(Z_M)  = \mathrm{Krull 
dim}(R_{\mathfrak{m}})  =  \mathrm{dim}_{k}(\mathfrak{m}/\mathfrak{m}^2)
 = \mathrm{dim}_{k}(M_{\mathfrak{m}}/M^{2}_{\mathfrak{m}}),
\]
where the second equality holds since $R_{\mathfrak{m}}$ is regular, 
\cite[Section 14, p.104]{mat}, and the third since $M$ is unramified. Thus $Z_M$ is regular, as required.
\end{proof}    
\subsection{The Azumaya locus}
\label{azloc}
In this paragraph and indeed for the rest of this paper we shall consider the 
following set-up. We continue with the embedding of affine commutative 
$k$-algebras $\theta:R\longrightarrow Z$ of the previous paragraph. In addition 
we now assume that $Z$ is the centre of a prime $k$-algebra $T$ which is a 
finite $Z$-module, (so that $Z$ is a domain). We 
define the \textit{Azumaya locus} to be the set of 
maximal ideals $M$ of $Z$ such that $MT$ is a maximal ideal of $T$. Recall 
\cite[Proposition 3.1]{brogoo} that the simple $T$-modules   
all have finite $k$-dimension and there is a bound on these dimensions, 
namely the PI-degree of $T$; this bound is always attained. By Schur's lemma, 
$Z$ acts on simple $T$-modules by scalars, so to each simple 
$T$-module $V$ we can associate the central character $\zeta_V:Z 
\longrightarrow k$ defined by these scalars.

The  \textit{Azumaya locus of T} is the set
\[
\mathcal{A}_T = \{ \mathrm{ker}(\zeta_V) : V \textit{ a simple $T$-module of 
maximal 
dimension} \} ,
\]
a non-empty open subset of $\mathrm{Maxspec}(Z)$. We also define \textit{the fully Azumaya locus of $R$ with respect to $T$} to be
\[
\mathcal{F}_T = \{\mathfrak{m}\in \mathrm{Maxspec}(R) \quad : \quad \theta^{*^{-1}}(\mathfrak{m}) \subseteq \mathcal{A}_T \}.
\]
Using the facts that $\mathcal{A}_{T}$ is open and that images of closed sets under finite morphisms are closed \cite[Ex.II.3.5(b)]{har}, it is easy to show that $\mathcal{F}_{T}$ is a non-empty open subset of $\mathrm{Maxspec}(R)$. 
Let $V$ be a simple $T$-module with $\mathrm{ker}(\zeta_V)\in\mathcal{A}_T$. Notice that, 
since $\mathrm{ker}(\zeta_V)T$ is maximal, if $V'$ is another simple $T$-module  
then $V\cong V'$ if and only if $\zeta_V =\zeta_{V'}$. 
Secondly, suppose $\x{T}{1}{V}{V'}\neq 0$ and let $M=\ker (\zeta_V)T$ and 
$M'=\ker (\zeta_{V'})T$. Then $\x{T}{1}{T/M}{T/M'}\neq 0$, and since this $Z$-module 
is annihilated by the two maximal ideals $M \cap Z$ and $M' \cap Z$  we must 
have $M=M'$ and so $V\cong V'$. That is, $\x{T}{1}{V}{V'}\neq 
0$ only if $V\cong V'$. Of course the reverse implication is also true unless 
$Z_{\mathfrak{m}} = 0$. Similar remarks apply to $\x{T}{1}{V'}{V}$, and to 
higher Ext-groups.

\begin{prop}
\label{genaz}
Let $M$ be a maximal ideal of $Z$ and let $\mathfrak{m}=M\cap R$. Suppose that 
$M \in \mathcal{A}_T $. Let $\mathfrak{S} = 
Z\setminus M$, an Ore set in $T$, and set $T_M = \mathfrak{S}^{-1}T$. Let $n$ be 
the PI degree of $T$. There is an 
algebra isomorphism
\[
\frac{T_M}{\mathfrak{m}T_M} \cong \mathrm{Mat}_{n}\left( 
\frac{Z_M}{\mathfrak{m}Z_M} 
\right),
\]
where $T_M/\mathfrak{m}T_M$ is a direct summand of $T/\mathfrak{m}T$ and the right side of the isomorphism consists of $n\times n$ matrices over the $(M/\mathfrak{m}Z)$-primary component of $Z/\mathfrak{m}Z$.
\end{prop}
\begin{proof}
First note that since $M$ is Azumaya $T_M$ is a free 
$Z_{M}$-module, \cite[1.8.31]{row}. It follows that $Z_M/\mathfrak{m}Z_M$ is a 
(central) subalgebra of $T_M/\mathfrak{m}T_M$.

By definition $MT_M$ is a maximal ideal of $T_M$ and we have an algebra 
isomorphism
\begin{equation}
\label{eq:weemat}
\frac{T_M}{MT_M} \cong \mathrm{Mat}_{n}(k).
\end{equation}
As the image of $MT_M$ in 
$T_M/\mathfrak{m}T_M$ is nilpotent we see that $T_M/\mathfrak{m}T_M$ is a local, 
finite dimensional algebra with simple right $T_M/\mathfrak{m}T_M$-module $V$ 
having dimension $n$.

Let $P(V)$ be the $T_M/\mathfrak{m}T_M$-projective cover of $V$. Then, as right 
$T_M/\mathfrak{m}T_M$-modules, we have
\begin{equation}
\label{eq:projsum}
\frac{T_M}{\mathfrak{m}T_M} \cong \bigoplus^{n} P(V).
\end{equation}
Therefore there are algebra isomorphisms
\[
\frac{T_M}{\mathfrak{m}T_M}\cong 
\mathrm{End}_{\frac{T_M}{\mathfrak{m}T_M}}\left( 
\bigoplus^{n}P(V)\right) \cong 
\mathrm{Mat}_{n}\left(\mathrm{End}_{\frac{T_M}{\mathfrak{m}T_M}}\left(P(V)\right
)\right).
\]
Since $Z_M/\mathfrak{m}Z_M$ is central in $T_M/\mathfrak{m}T_M$ we have a map
\begin{equation}
\label{eq:emcen}
\frac{Z_M}{\mathfrak{m}Z_M} \longrightarrow 
\mathrm{End}_{\frac{T_M}{\mathfrak{m}T_M}}\left( P(V) \right),
\end{equation}
given by $z\longmapsto (p\longmapsto p.z)$. By (\ref{eq:projsum}) and the 
freeness of $T_M/\mathfrak{m}T_M$ as a $Z_M/\mathfrak{m}Z_M$-module, this map is 
injective. By (\ref{eq:weemat}) we see that as an $Z_M$-module, and so as an 
$R_{\mathfrak{m}}$-module, there is an isomorphism $T_M\cong 
\bigoplus^{n^2}Z_M$. Therefore
\[
\dim \left(\frac{T_M}{\mathfrak{m}T_M}\right) = 
n^2\dim\left(\frac{Z_M}{\mathfrak{m}Z_M}\right),
\]
so (\ref{eq:emcen}) is an isomorphism.
\end{proof}

\subsection{}
We deduce the following.
\begin{cor}
Let $M$ be a maximal ideal of $Z$ lying on the Azumaya locus and let 
$\mathfrak{m}=R\cap M$. Then the ramification index of $M$ equals the Loewy 
length of the finite dimensional algebra $T_M/\mathfrak{m}T_M$.
\end{cor}
\subsection{}
\label{regalgs}
The following  special case of Proposition \ref{genaz} will prove useful later. It follows from the proposition because under the stated hypothesis $T/\mathfrak{m}T$ is the direct sum of its primary components, and the latter are $n\times n$ matrix rings over the primary components of $Z/\mathfrak{m}Z$. Recall from (\ref{azloc})the definition of the fully Azumaya points, $\mathcal{F}_T$.
\begin{cor}
Suppose $\mathfrak{m}\in \mathcal{F}_T$. Then there 
is an 
isomorphism
$$
\frac{T}{\mathfrak{m}T}\quad \cong \quad \mathrm{Mat}_n\left( 
\frac{Z}{\mathfrak{m}Z}\right).
$$
\end{cor} 
\subsection{}
\label{azunram}
In general the inclusion given by Lemma \ref{smooth} is strict. The 
following result makes this inclusion precise, in representation-theoretic 
terms, in the presence of an additional hypothesis whose validity in various 
important cases is discussed in (\ref{az=sing}).  
\begin{thm}
Suppose that $\mathrm{Maxspec}(R)$ is smooth and that the smooth locus of 
$\mathrm{Maxspec}(Z)$ coincides with the Azumaya locus of $T$. Let 
$M\in\mathrm{Maxspec}(Z)$ and let $\mathfrak{m}=R\cap M$. Then the following are 
equivalent
\newline
(i) $M$ is unramified;
\newline
(ii) $M$ is on the Azumaya locus and $T/MT$ is a projective 
$T/\mathfrak{m}T$-module.
\end{thm}
\begin{proof}
Let $\mathfrak{S}=Z\setminus M$, a multiplicatively closed Ore set. For every 
$T/MT$-module $N$, there is a 
natural equivalence
\[
\mathfrak{S}^{-1}\hm{T/\mathfrak{m}T}{N}{-} \cong 
\hm{\mathfrak{S}^{-1}(T/\mathfrak{m}T)}{\mathfrak{S}^{-1}N}{\mathfrak{S}^{-1}-}. 
\]
Multiplication by elements of $\mathfrak{S}$ on the $T/MT$-module $N$ 
is bijective, so the above isomorphism yields
\[
\hm{T/\mathfrak{m}T}{N}{-}\cong \hm{\mathfrak{S}^{-1}(T/\mathfrak{m}T)} 
{\mathfrak{S}^{-1}N}{\mathfrak{S}^{-1}-}.
\]
Combining this with Proposition \ref{genaz} proves that if $M$ is on the Azumaya 
locus then $T/MT$ is a projective $T/\mathfrak{m}T$-module if and only if 
$\mathfrak{S}^{-1}(T/MT)=T_M/MT_M$ is a projective $Z_M/\mathfrak{m}Z_M$-module.
It's clear from Proposition \ref{genaz} that this happens if and only if   
$Z_M/\mathfrak{m}Z_M = k$, in other words if and only if $\mathfrak{m}Z_M = 
MZ_M$, that is $M$ is unramified.

Now by hypothesis the Azumaya locus of $\text{Maxspec}(S)$ coincides with the 
smooth locus of $\text{Maxspec}(S)$. Thus combining the above paragraph with 
Lemma \ref{smooth} yields the theorem.
\end{proof}

\subsection{}
\label{az=sing}
The hypothesis of Theorem \ref{azunram} that the smooth and Azumaya loci 
coincide is rather restrictive - for example it fails to hold for the 
2-dimensional solvable non-Abelian Lie algebra in positive characteristic, and 
for $U(\mathfrak{s}\mathfrak{l}(2))$ in characteristic 2, \cite[Examples 
3.4]{brogoo}. However it is 
satisfied in three important cases which we list in the theorem below, the first two of which 
we shall explore in detail in the remaining sections of the paper. Full details 
of the algebras listed below will be given for cases 1 and 2 in
paragraphs (\ref{liedefs}) and (\ref{notqg}).  Here and throughout
this paper we shall assume in case 1 that $\g$ is the Lie algebra of a
connected, reductive group $G$ over an algebraically closed field of
characteristic $p$ which satisfies the following hypotheses:
\begin{list}{}{} 
\item (A) the derived group $\mathcal{D}G$ of $G$ is simply-connected;
\item (B) $p$ is odd and a good prime for $G$;
\item (C) the trace form on $\g$ is non-degenerate.
\end{list}
We remark that, in the first part of
the following theorem, only the  
case $\g$ semisimple, $p$ very good is considered in \cite{brogoo}. We
shall explain in Section 3 how to weaken the hypotheses to those stated here. Hypotheses (B) and (C) are discussed in \cite[6.3,6.4,6.5]{jan6}. The prime $p$ is good for $\g$ if and only if it is good for all the irreducible components of the root system of $\g$. The primes which
are not good for an irreducible root system are
\begin{itemize}
\item $p=2$ for types $B_r$, $C_r$ and $D_r$;
\item $p=2$ or $3$ for types $E_6$, $E_7$, $F_4$ and $G_2$;
\item $p=2,3$ or $5$ for type $E_8$.
\end{itemize}
Hypothesis (C) entails the exclusion of $\g$ simple of type $A_r$ when $p|r+1$. Crucially, however, both hypotheses are preserved under passage to Levi factors.

For the function algebras
$\EO$ the questions addressed for cases 1 and 2 in the present paper
are considered in \cite{brogor2}.
\begin{thm}\cite{brogoo}
The hypotheses of Theorem \ref{azunram} are satisfied in the following cases:
\newline
1. Suppose that hypotheses (A), (B) , and (C) listed above hold. Take $T= U(\g)$, $Z=Z(U(\g))$, and $R=Z_0$, the 
$p$-centre of $U(\g)$.
\newline
2 Let $\g$ be a complex semisimple Lie algebra, and let $\epsilon$ be an $\ell 
th$ root of unity, where $\ell$ is odd and prime to 3 if $\g$ involves a factor 
of type $G_2$. Take $T$ to be the simply connected quantised enveloping algebra $U_{\ep}(\g)$, 
$Z=Z(U_{\ep}(\g))$, $R = Z_0(\g)$.
\newline
3. Let $\g$ and $\ell$ be as in 2 and let $G$ be the connected, simply connected 
algebraic group with Lie algebra $\g$. Take $T=\EO$, $Z=Z(\EO)$, $R = \mathcal{O}[G].$
\end{thm}

\subsection{Blocks}
\label{Blocks}
We continue with a triple $R \subseteq Z \subseteq T$ as introduced in
(\ref{azloc}). Let $\mathfrak{m}$ be a maximal ideal of $R$. A basic
problem in the representation theory of $T$ is to describe the blocks
of the finite dimensional algebra $T/\mathfrak{m}T$, which are in
one-to-one correspondence with the primitive central idempotents of
the latter algebra. A particularly convenient scenario for determining
these idempotents occurs when the centre of $T/\mathfrak{m}T$ equals
$Z + \mathfrak{m}T/\mathfrak{m}T$. Unfortunately this is not always
the case, and in fact we include a counterexample in the enveloping
algebra setting in (\ref{Muller}). However it turns out that the
failure of this equality is always confined to the nilpotent elements
of the centre of $T/\mathfrak{m}T$. This is the content in our setting
of the following result of B. M\"{u}ller, \cite[Theorem 7]{mul}, \cite[Theorem 11.20]{goo-war}, which is sufficient to enable us to get good information on the blocks of $T/\mathfrak{m}T$ (which of course correspond precisely as $\mathfrak{m}$ varies to the blocks of $T$, since $\mathfrak{m}T$ is centrally generated).
\begin{prop}
\label{Muller}
(M\"{u}ller) Let $R$, $Z$ and $T$ be as in (\ref{azloc}) (although there is no need here to assume that $R$ is regular), and let $\mathfrak{m}$ be a maximal ideal of $R$. The primitive central idempotents of $T/\mathfrak{m}T$ are the images of the primitive idempotents of $Z/\mathfrak{m}Z$. 
\end{prop}
\begin{cor}
Let $R,\quad Z,\quad T$ and $\mathfrak{m}$ be as in the proposition. Then the blocks of $T/\mathfrak{m}T$ are in one-to-one correspondence with the maximal ideals of $Z$ lying over $\mathfrak{m}$.
\end{cor} 
\section{Enveloping algebras in positive characteristic}
\subsection{}
\label{liedefs}
Let $G$ be a connected, reductive algebraic group over $k$, an algebraically
closed field of characteristic $p$. We shall continue throughout
Section 3 to assume that $\g=\text{Lie}(G)$ and $p$ satisfy hypotheses
(A), (B) and (C) of (\ref{az=sing}). Let $T$ be a maximal torus of $G$
and let $\h = \mathrm{Lie} (T)$. Let $\Phi$ be the root system of
$G$ with respect to $T$. For each $\alpha\in \Phi$ let $U_{\alpha}$ denote the
corresponding root subgroup of $G$ and let
$\g_{\alpha}=\text{Lie}(U_{\alpha})$ be its Lie algebra, a root
subspace of $\g$. We will abuse notation by considering
$\alpha \in \h^*$ rather than its proper designation $d\alpha$.
 Choose a system $\Phi^+$ of positive roots and set $\n^+$
equal to the sum of all $\g_{\alpha}$ with $\alpha > 0$. For $\alpha \in \Phi^+$, let $h_{\alpha}\in [\g_{\alpha},\g_{-\alpha}]$ be the
unique element of $\h$ such that $\alpha (h_{\alpha})=2$ (recall $p$ is odd). The
subalgebra $\n^-$ is similarly defined on $\Phi^-$, the negative
roots. We have the triangular decomposition
\[
\g =  \n^- \oplus \h \oplus \n^+.
\]
Let $\mathfrak{b}^+ = \h \oplus \n^+$, the Lie algebra of a Borel subgroup of $G$
containing $T$. Let $\Delta = \{\alpha_1, \ldots ,\alpha_r \}$ denote the simple roots associated with
the choice of positive roots $\Phi^+$.

For each root $\alpha$ fix a basis vector
$x_{\alpha}\in \g_{\alpha}$. Then $x_{\alpha}^{[p]}=0$ for all $\alpha
\in \Phi$. Since $T$ is a torus $\h$ is an abelian Lie algebra and has
a basis $\{ h_1,\ldots ,h_r\}$ such that $h_i^{[p]} = h_i$. In
particular the rank of $G$ is $r$, and if we set $2N = |\Phi|$ then
$\g$ has dimension $2N+r$.

Let $X=X(T)$ be the character group of $T$. This contains the root
lattice, 
$Q=\mathbb{Z}\Phi$, as a subgroup. Let $\Lambda$ be the character group
reduced modulo $p$, that is $X/pX$. There is an inclusion of $X/pX$
into $\h^*$,
where $X/pX$ is identified with set $\{ \lambda \in \h^* : \lambda
(h_i)\in \mathbb{F}_p\}$, \cite[11.1]{jan6}. Since the $h_{\alpha}$
with $\alpha$ simple are linearly independent in $\h$ by Hypothesis
(A) we can find $\rho\in
\h^*$ with $\rho (h_{\alpha})=1$ for all $\alpha\in \Delta$. We fix
such a $\rho$ once and for all. 

We write $gx$ for the adjoint action of an element $g\in G$ on an
element $x\in \g$. Similarly we will write $g\chi$ for the coadjoint
action of $G$ on $\g^*$, defined by $g\chi (x) = \chi
(g^{-1}x)$. Let $W$ be the Weyl group of $G$. Then $W$ is generated by
the simple reflections $s_{\alpha}$ for all $\alpha\in \Delta$. There
is an action of $W$ on $\h^*$ given by $s_{\alpha}(\lambda) = \lambda
- \lambda(h_{\alpha}) \alpha$ for all $\lambda\in \h^*$ and $\alpha\in
\Delta$.

By considering weights we see that $\g_{\alpha}$ is orthogonal to each
$\g_{\beta}$ with $\beta\neq -\alpha$ with respect to the non-degenerate form of
Hypothesis (C). Therefore this form induces a $W$-invariant non-degenerate form on $\h$
by restriction. 
\subsection{The centre of $U(\g)$.}
\label{centresred}
Let $U = U(\g)$ be the enveloping algebra of $\g$ and let $Z=Z(\g)$ be the centre of 
$U(\g)$. The main features of the description of $Z$ go back to the
work of Veldkamp \cite{veld} in the 70s, improvements to which have
been made in \cite{mirrum}, who obtained results for $\g$ semisimple
and $p$ very good. However, some of the results
of \cite{veld} do not extend to our setting. Indeed if $\g = \mathfrak{sl}(p+1)$ then the
subalgebra of $\g$ consisting of all elements which commute with
$e_{11}+e_{22}+\ldots +e_{pp}\in\g$ provides a counterexample to the natural extension of
\cite[Theorem 6.3]{veld}. For all Lie algebras satisfying Hypotheses
(A), (B) and (C), however, we will show that the description of $Z$
found in \cite[Theorem 3.1]{veld} and \cite{mirrum} remains valid. To
prove this requires preparation.

Recall there is a Jordan decomposition in $\g$: each element $x\in \g$
can be written uniquely as $x=x_s + x_n$ with $x_s$ semisimple, $x_n$
nilpotent and $[x_s,x_n] = 0$. Given $x\in \g$ we can always find $g\in
G$ such that $g.x \in \bp$, \cite[Proposition 14.25]{bor}. Let $\mathfrak{c}_{\g}(x) = \{ y \in \g
: [x, y]= 0\}$ and $C_G(x) = \{ g\in G: g.x = x\}$. An element $x \in \g$ is called \textit{regular} if
$\dim (C_G(x)) = r$, the rank of $G$.

Let $\theta :\g \longrightarrow
\g^*$ be the $G$-invariant isomorphism induced by the non-degenerate
form of Hypothesis (C). We use this to transfer the Jordan decomposition from $\g$ to
$\g^*$. In particular any element of $\g^*$ is conjugate to $\chi\in
\g^*$ such that $\chi (\mathfrak{n^+}) = 0$. For $\chi \in \g^*$ let $\cen{\chi} = \{ y\in \g
: \chi ([\g , y]) = 0\}$ and $C_G(\chi) = \{ g\in G: g.\chi =
g\}$. It is easy to check that $C_G(x) = C_G(\theta (x))$ and
$\cen{x} = \cen{\theta(x)}$.

\begin{lem}
Let $x\in \g$. Then $\text{Lie}(C_G(x)) = \cen{x}$. Moreover, if $x\in
\h$ then $C_G(x)$ is a connected, reductive algebraic group of rank
$r$ satisfying
hypotheses (A), (B) and (C). Indeed, in this case $C_G(x)$ is generated by $T$ and
the root subgroups $U_{\alpha}$ with $\alpha (x) = 0$.
\end{lem}
\begin{proof}
The existence of the
non-degenerate bilinear form on $\g$ allows us to identify
$\text{Lie}(C_G(x))$ and $\cen{x}$ by \cite[Theorem
3.10]{humCC}. That $C_G(x)$ is a connected, reductive algebraic group follows from
\cite[II.3.19]{S-S} and our assumptions (A) and (B). The precise
description of $C_G(\chi)$ follows from \cite[7.4]{jan6}.  This description shows that it is a
Levi subgroup and so, by \cite[6.5]{jan6}, satisfies the hypotheses. 
\end{proof}
\indent Standard arguments, \cite[Chapter 4]{humCC}, \cite[Theorem 3.3(a)]{S-S}, show that regular nilpotent elements exist in
$\g$ and form a single class under the adjoint action of $G$. Moreover
Hypothesis (C) ensures there are only a finite number of nilpotent
classes in $\g$, \cite[Corollary 4.2]{ric}, \cite[Theorem 3.10]{humCC}. With this in hand the results of \cite[Section
4]{veld} follow verbatim, giving us the following proposition.
\begin{prop}
The regular elements of $\g$ form a Zariski open set with complement
of codimension at least 3.
\end{prop}
\begin{rems}
(i) It can be shown that the regular semisimple elements of $\g$ form a non-empty Zariski
open set of $\g$; see \cite[Theorem 2.5]{humCC} for the group case.

(ii) Let $x= x_s + x_n\in \g$. It is clear that
$\cen{x}=\mathfrak{c}_{\cen{x_s}}(x_n)$. Therefore, by Lemma \ref{centresred},
$x$ is regular in $\g$ if and only if $x_n$ is regular in
$\cen{x_s}$. In particular this shows that given any semisimple
element $y\in \g$ there exists a regular element $x\in \g$ such that
$x_s = y$.
\end{rems}
\subsection{}
\label{Ginv}
We now turn to the adjoint quotient $\pi :\g \longrightarrow \g/G$,
the morphism associated with the inclusion $S(\g)^G \longrightarrow
S(\g)$. Our objective is Corollary \ref{adjquot}, which is the key tool in deriving the main properties of the centre listed in Theorem \ref{cent}. The proof we offer of Corollary \ref{adjquot} closely follows the argument given by Chriss and Ginzburg \cite[Section 6.7]{C-G} to prove the corresponding result of Kostant in characteristic zero. We shall therefore simply refer the reader to \cite{C-G} at points where the argument is identical. 
\begin{lem}
The restriction map from $S(\g)$ to $S(\h)$ yields an isomorphism

\[
\psi : S(\g)^G\longrightarrow S(\h)^W.
\]
\end{lem}
\begin{proof}
The map
\[
\psi^* : Sym ( \g)^G \longrightarrow Sym (\h)^W,
\]
induced by orthogonal projection from $\g$ to $\h$, is an
isomorphism. To see this note that by \cite[Corollary 4.5]{veld} $f(\eta) = f(\eta_s)$ for any $f \in Sym
(\g)^G$ and $\eta \in \g^*$, with Jordan
decomposition $\eta = \eta_s + \eta_n$. Since $\eta_s$ is conjugate to element of $\h^*$ it
follows that $\psi$ is injective. Replacing $U(\g)$ by $S(\g)$ in
\cite[9.6]{jan6} shows that $\psi$ is indeed an isomorphism.  The lemma
follows by applying the non-degenerate form provided by Hypothesis (C).
\end{proof}
By \cite[Th\'{e}or\`{e}me 2 and Corollaire]{dem} and \cite[9.6]{jan6} there are homogeneous
invariants $T_1,\ldots ,T_r$ such that the algebra $S(\h)^W$ is a
polynomial algebra on these generators. Thus the lemma shows that
$S(\g)^G$ is generated by $J_1, \ldots J_r\in S(\g)^G$, the
preimages along $\psi$ of the elements $T_1,\ldots ,T_r$.  By \cite[Theorem
3.9]{humCOX} we have $\sum_{i=1}^r deg(J_i) = N + r$ where, as always, $N$ is the
number of positive roots of $G$. 

\subsection{}
\label{adjquot}
Coadjoint orbits in $\g^*$ occur with a natural symplectic
structure, \cite[Chapter 15, Theorem 1]{kir}. Indeed, let $\mathcal{O} = G.\chi$ be a coadjoint
orbit. Thanks to Lemma \ref{centresred} the isomorphism $G/C_G(\chi) \longrightarrow \mathcal{O}$ induces an
isomorphism between $\g/\cen{\chi}$ and $T_{\chi}\mathcal{O}$, see also 
\cite[3.8, 3.10]{humCC}. The required symplectic structure is
obtained from the skew symmetric form
\[
\omega_{\chi} : \g/\cen{\chi} \times \g/\cen{\chi} \longrightarrow k,
\]
which sends $(x,y)$ to $\chi([x,y])$. 

Through $\theta$, the $G$-invariant isomorphism from $\g$ to $\g^*$,
we can transfer our non-degenerate form of Hypothesis (C) from $\g$ to
$\g^*$. This induces a $G$-invariant non-degenerate form on the spaces $\Lambda^i\g^*$ \cite[1.6.5]{Darling}. Moreover,
since we are assuming $k$ does not have characteristic two, we can find an orthonormal
basis for $\g^*$ with respect to the form, \cite[1.6.3]{Darling}. Let $\{
X_1, \ldots ,X_n\}$ be such a basis. Of course, exterior products of
elements from this basis will yield orthonormal bases for the spaces
$\Lambda^i\g^*$, \cite[1.6.5.1]{Darling}.

Since $\dim (\Lambda^{n}\g^*)= 1$ there is a unique, up to sign,
$n$-form, say $V$, satisfying $(V, V) =1$. This can be realised with
the above basis as the element
$X_1\wedge \ldots \wedge X_n$. Consequently we can define the Hodge $\ast$-operator, a linear isomorphism 
\[
\ast : \Lambda^i\g^* \longrightarrow \Lambda^{n-i}\g^*,
\]
which is determined by requiring that, for all $\alpha, \beta \in \Lambda^i\g^*$,
\[
(\alpha , \beta ) V = \alpha \wedge (\ast \beta).
\]
Let $\Omega^{\scriptscriptstyle \bullet}(\g) = S(\g) \otimes
\Lambda^{\scriptscriptstyle \bullet}\g^*$ be the $S(\g)$-module
of polynomial differential forms on $\g$. Extending $\ast$ by
$S(\g)$-linearity gives an $S(\g)$-module homomorphism $\ast :
\Omega^i (\g) \longrightarrow \Omega^{n-i} (\g)$. We will write $dX_i$
when considering $X_i$ as a $1$-form in $\Omega^1(\g)$. 

For $x\in\g$ let $\Omega_x \in \Lambda^2\g^*$ be defined by $\Omega_x(y,z) = (x,
[y,z])$, for all $y,z\in \g$. The assignment $\Omega : x \longmapsto
\Omega_x$ gives an element of $\Omega^2(\g)$, whose polynomial coefficients with
respect to the basis $\{X_1, \ldots X_n\}$ are linear. Let $\Omega^N =
\Omega \wedge \ldots \wedge \Omega$ be the $N^{\text{th}}$ exterior
power of $\Omega$.

\begin{lem}
A point $x\in \g$ is regular if and only $\Omega_x^N \neq 0$.
\end{lem}
\begin{proof}
For $x\in \g$ let $\chi = \theta (x)$ and let $\mathcal{O}=G.\chi$. As observed
above there is an isomorphism between $T_{\chi}\mathcal{O}$ and $\g /
\cen{x}$. By definition $\Omega_x(y,z) = (x, [y,z]) = \chi ([y,z]) =
\eta^*\omega_{\chi} (y,z)$ for all $y,z\in\g$, where $\eta : \g
\longrightarrow \g/\cen{\chi}$ is the natural map. Therefore $\Omega_x^N \neq 0$ if and only if
$\bigwedge^N\omega_{\chi}\neq 0$. But since $\omega_{\chi}$ is
non-degenerate on $T_{\chi}\mathcal{O}$ and $\dim \mathcal{O}\leq 2N$
this happens if and only if $\dim \mathcal{O} = 2N$, in other words
$\chi$ (and hence $x$) is regular.
\end{proof}
\indent We shall need to use the final part of the following lemma in the proof of Proposition \ref{adjquot}. Since we could find no suitable reference and the result is wrongly stated in \cite{C-G}, we include here some indications of the proofs. Recall that if $W$ is an $k$-vector space and $x \in W, f \in \wedge^pW^*$, the \textit{interior product} $x \rightharpoonup f$ is $f(x \wedge -) \in \wedge^{p-1}W^*$, (or $f(x)$ if $p=1$). We shall write $\mathrm{rad}(f) = \{x \in W : x \rightharpoonup f = 0\}$. As above, the form on $W$ induces non-degenerate forms (which we denote by the same symbol $(\quad,\quad)$) on $W^*$ and its exterior powers. We use the symbol $\perp$ to denote orthogonality with respect to any of these induced forms. 
\begin{lem}
Let $W$ be an $n$-dimensional $F$-vector space admitting a
non-degenerate bilinear form $( \quad,\quad )$ (recall that the
characteristic of $k$ is not two). Fix $1 \leq p,q \leq n$, and let $x \in W, f \in \wedge^{p}W^*$ and $g \in \wedge^qW^*$. Let $\{v_1, \ldots , v_n\}$ be an orthonormal basis for $W$ with dual basis $\{v_{1}^*, \ldots , v_{n}^*\}$ and isomorphism $\theta:v_i \longrightarrow v_i^*$. For $I = \{i_1, \ldots , i_p : i_1 < i_2 < \ldots < i_p \}$ write $f_I$ to denote the element $v_{i_1}^* \wedge \ldots \wedge v_{i_p}^*$ of $\wedge^{p}W^*$.\\
1. $\ast f_I = (-1)^{sgn(I|J)}f_J$, where $J = \{1, \ldots ,n\} - I$ and $sgn(I|J)$ denotes the sign of the permutation sending 1 to $i_1, \ldots , n$ to $j_{n-p}$.\\
2. $\ast(x \rightharpoonup f) = (-1)^{p+1}\theta(x)\wedge(\ast f)$.\\
3. $ x \rightharpoonup (f \wedge g) = (x \rightharpoonup f) \wedge g  +  (-1)^pf \wedge (x \rightharpoonup g) $.\\
4. $\mathrm{rad}(\ast f) \subseteq \mathrm{rad}(f)^{\perp}$.
\end{lem}
\begin{proof}
1. Straightforward from the definitions.\\
2. Prove this first for $x = v_i$ and $f = f_I$. Both sides are zero if $i$ is not in $I$. If $i = i_s \in I$, set $I_s = I - \{i_s\}$ and prove $(v_{i_s} \rightharpoonup f_I) = (-1)^{s+1}f_{I_s}$. Then apply part 1. The general case follows by bilinearity.\\
3. By linearity one need only treat the case where $f$ and $g$ are exterior products of $1$-forms. The identity follows for such $f$ and $g$ by expanding $x \rightharpoonup(f \wedge g)$ as a determinant down the first column.\\
4. Suppose that $x \in \mathrm{rad}(f), y \in \mathrm{rad}(\ast f)$. Then $\theta(x) \wedge (\ast f) = 0$ by part 2. Thus $0 = y \rightharpoonup(\theta(x) \wedge (\ast f)) = (x,y)(\ast f)$, the second equality using part 3. That is, $x \in (\mathrm{rad}(\ast f))^\perp$, so $\mathrm{rad(f)} \subseteq (\mathrm{rad}(\ast f))^\perp$. Applying this inclusion to $\ast f$ and using part 1 yields the desired conclusion.
\end{proof}  
 
Recall that in $\Omega^1(\g)$ we have
\[
dJ_i = \sum_{j=1}^n \frac{\partial J_i}{\partial X_j}dX_j.
\]
\begin{prop}
We have
\begin{equation}
\label{eq:regg}
\ast (\Omega^N) = C\cdot dJ_1\wedge \ldots \wedge dJ_r
\end{equation}
where $C$ is a non-zero constant.
\end{prop}
\indent The proof follows exactly the argument for the corresponding statement in characteristic zero given in \cite[Theorem 6.7.32]{C-G}. At one point Chriss and Ginzburg appeal to the identity $\mathrm{rad}(\ast f) = (\mathrm{rad}(f))^\perp$, which is false in all characteristics. However only the inclusion provided by part 4 of the above lemma is needed. 
\begin{cor}
The adjoint quotient map
\[
\pi : \g \longrightarrow \g/G
\]
is smooth on the regular elements of $\g$, a Zariski open set with
complement of
codimension at least three.
\end{cor}
\begin{proof}
Under the identification of $S(\g)^G$ with a polynomial algebra in $r$
variables given by Lemma \ref{Ginv} the map $\pi$ can be written as 
\[
\pi(x) = (J_1(x),\ldots ,J_r(x)),
\]
for all $x\in\g$. Since $\g$ and $\g/G$ are smooth varieties $\pi$ is
smooth at $x$ if and only if the differential of $\pi$ at $x$ is
surjective, \cite[Chapter III, Proposition 10.4]{har}. In the above
coordinates the differential at $x$ is described by the matrix
\[
\left( \frac{\partial J_i}{\partial X_j} \right)_{\substack{1\leq
i\leq r \\ 1\leq j\leq n}}(x).
\]
The first of the above lemmas and the proposition show that this matrix has rank
$r$ if and only if $x$ is a regular element, proving the first
claim. The second part is Proposition \ref{centresred}.
\end{proof} 
\subsection{}
\label{cent}
Let $Z_0 = k[x^p - x^{[p]} : x\in \g]\subseteq Z$ be the $p$-centre 
of $U$, a polynomial algebra in $\dim(\g)$ indeterminates. It follows from 
the Poincare-Birkhoff-Witt theorem that $U(\g)$ is a 
free $Z_0$-module of rank $p^{\dim (\g)}$. As in section (\ref{nonemp}), we denote the
embedding of $Z_0$ in $Z$ by $\theta$. Let $Z_{1}$ be the ring of
invariants of $U(\g)$ under the adjoint action of $G$, so $Z_{1}$ is
also clearly a subalgebra of $U(\g)$.

For a $k$-vector space $V$ let $V^{(1)}$ denote the Frobenius twist of $V$, that 
is the 
$k$-vector space with the same elements as $V$ but with $k$ acting via 
$a\longmapsto 
a^{1/p}$. Then there is an isomorphism $\Psi:Z_0 \longrightarrow
\mathcal{O}(\g^{*(1)})$ defined by $\Psi(x^p-x^{[p]})(\eta)=\eta(x)$
for $x\in \g$ and $\eta \in \g^{*(1)}$. Under this isomorphism
$Z_0\cap Z_1=Z_0^G$ is sent to $\mathcal{O}(\g^{*(1)})^G$, so the
inclusion $Z_0\cap Z_1 \longrightarrow Z_0$ corresponds to the
coadjoint quotient $\g^{*(1)}\longrightarrow \g^{*(1)}/G$. 
\begin{thm}
\label{centreU}
Retain the notation introduced in (\ref{liedefs}) and above, so $\g$ is a reductive Lie algebra. Assume that hypotheses (A), (B) and (C) of (\ref{az=sing}) hold.\\
1.  Let $\gamma: S(\h)\longrightarrow S(\h)$ be defined
on generators as $\gamma(h) = h - \rho (h)$. Then the Harish-Chandra map
\begin{equation}
\label{eq:HCclass}
\Theta := \gamma \circ(\ep \otimes Id_{U(\h)}\otimes \ep):U\cong U(\n^-)\otimes
U(\h)\otimes U(\n^+) \longrightarrow U(\h)\cong S(\h),
\end{equation}
restricts to an algebra isomorphism from $Z_1$ to $S(\h)^W$. Recall
that $T_{1}, \ldots , T_{r}$ are homogeneous invariants such that $S(\h)^{W}$ is the polynomial algebra on these generators.\\
2. Write $Y_{i} = \Theta^{-1}(T_{i})$ for $i = 1, \ldots , r$. Then $Z_{0} \cap Z_{1}$  is a polynomial algebra of rank $r$, and $Z_{1}$ is a free $Z_{0} \cap Z_{1}$-module with basis $\{Y_{1}^{m_{1}} \ldots Y_{r}^{m_{r}} : 0 \leq m_{i} < p, 1 \leq i \leq r \}$.\\
3. $Z_{1}$ is a complete intersection over $Z_{0} \cap Z_{1}$.\\
4. $Z_{0} \otimes_{Z_{0} \cap Z_{1}} Z_{1}$ is a complete intersection ring which is smooth in codimension 1.\\
5. The multiplication map $\mu :Z_{0} \otimes_{Z_{0} \cap Z_{1}} Z_{1}\longrightarrow Z$ is an isomorphism of algebras.\\
6. $Z$ is a free $Z_{0}$-module of rank $p^r$.    
\end{thm}
\begin{proof}
1. This follows from \cite[Theorem 9.3]{jan6}, since the required
hypotheses hold by the final sentence of \cite[9.6]{jan6}.

2 and 6. We apply \cite[Theorem 3.1]{veld}, which is valid here in
light of Lemma \ref{Ginv} and Corollary \ref{adjquot}. We conclude that $Z$
is generated by $Z_0$ and $Z_1$ and in fact that $Z$ is a free
$Z_0$-module with basis $\{ Y_1^{i_1}\ldots Y_r^{i_r}: 0\leq i_j \leq
p-1 \text{ for all }j\}$. This proves 6, whilst taking $G$-invariants proves
the final claim of 2. That $Z_0\cap Z_1=Z_0^G$ is a polynomial algebra
on $r$ generators follows from Lemma \ref{Ginv}.

3. We have a commutative diagram 
\[
\begin{CD}
Z_0^G=Z_0\cap Z_1 @> \Theta >> k[h^p-h^{[p]}: h\in \h]^W\\
@V \Psi VV @VV \Psi| V \\
\mathcal{O}(\g^{*(1)})^G @> \text{res} >> \mathcal{O}(\h^{*(1)})^W.
\end{CD}
\]
where we
have noted that $\gamma (h^p-h^{[p]}) = h^p - h^{[p]}$ for all $h\in \h$. By Lemma
\ref{Ginv} the bottom map is an isomorphism which shows that the whole diagram consists of
isomorphisms since $\Psi$ is also an isomorphism. So to prove that $Z_1$ is a complete intersection over
$Z_0\cap Z_1$ it is enough to show that $k[h:h\in \h]^W$ is a complete
intersection over $k[h^p-h^{[p]}:h\in \h]^W$. To see this note that
there is a $W$-invariant $k$-algebra isomorphism $\psi$ from $k[h^p: h\in \h]$ to
$k[h^p-h^{[p]}:h\in\h ]$  obtained by sending $h_i^p$ to
$h_i^p-h_i$. Now it's easy to check that  $k[h^p: h\in \h]^W = k[T_1^{p} , \ldots, T_r^{p}]$, and therefore $k[h^p-h^{[p]}:h\in \h ]^W =
k[T_1^{(p)} , \ldots, T_r^{(p)}]$ where $T_i^{(p)} = T_i(h_1^p-h_1 ,
\ldots ,h_r^p-h_r) = \psi (T_i)$. Noting that $T_i^{(p)} = T_i^p + q_i(T_1,\ldots
T_r)$, where the degree of $q_i$, considered as an element of $U(\h)$,
is strictly less than $p\text{deg}(T_i)$, a straightforward degree
argument shows that 
\[
k[h:h\in \h]^W \cong \frac{k[h^p-h^{[p]}:h\in \h][X_1, \ldots
X_r]}{(T_i^{(p)} - X_i^p - q_i(X_1,\ldots X_r):1\leq i\leq r)},
\]
as required.

4. and 5. These follow exactly as in \cite[Theorem 6.4]{deckacpro4} in light of
parts 2. and 3. and Corollary \ref{adjquot}.
\end{proof}

\subsection{Proof of Theorem \ref{az=sing}.1.}
\label{Proofaz=sing}
It follows immediately from Proposition \ref{centresred}, Premet's
theorem \cite[Theorem 3.10]{pre2}, \cite[Theorem 7.6]{jan6}, and \cite[Proposition 3.1]{brogoo} that
\begin{equation}
\label{nonaz}
\textit{the locus of non-Azumaya points has codimension at least 3 in }Z.
\end{equation}
Using (\ref{nonaz}) in place of \cite[Proposition 4.8(ii)]{brogoo} one can now proceed exactly as in \cite[(4.10)]{brogoo} to prove Theorem \ref{az=sing}.1. 
\subsection{}
\label{frob}

Recall from Section \ref{cent} the isomorphism $\Psi: Z_0 \longrightarrow
\mathcal{O}(\g^{*(1)})$. For $\chi \in \g^{*}$ define $\chi^p\in \g^{*(1)}$ by
$\chi^p(x)=\chi(x)^p$ for $x\in\g$. Let $\mathfrak{m}_{\chi}$ be the inverse image under $\Psi$ of the zeros of $\chi^p$ in $\mathcal{O}(\g^{*(1)})$. Thus $\mathfrak{m}_{\chi} = \langle x^p - x^{[p]} - \chi(x)^p : x \in \g \rangle$, and every maximal ideal of $Z_0$ has this form. Set $U_{\chi} = U(\g)/\mathfrak{m}_{\chi}U(\g)$, an algebra of dimension $p^{\mathrm{dim}(\g)}$. For details of all this, see for example \cite[(2.3)-(2.10)]{jan6}. As the isomorphism class of $U_{\chi}$ depends only on the
$G$-orbit of $\chi \in \g^{*}$, \cite[2.9]{jan6}, it suffices, by
Section \ref{centresred}, to look at $\chi$ satisfying
$\chi (\n^{+}) = 0$. In this case if $\chi = \chi_s+\chi_n$ is the
Jordan decomposition we also have
$\chi_s(\n^{-})=0$ and $\chi_n (\h)= 0$ and so in particular
we can consider $\chi_s \in \h^{*}$. 

Now let $\chi \in \g^{*}$ with $\chi(\mathfrak{n}^{+})= 0$, and set
\begin{equation}
\label{lambdachi}
\Lambda_{\chi} = \{\lambda \in  
\mathfrak{h}^* : \lambda (h)^p-\lambda(h^{[p]}) = 
\chi(h)^p, \textrm{ for all } h \in \h \},
\end{equation}
so $\Lambda_{\chi} = \lambda + \Lambda$ for any $\lambda \in \Lambda_{\chi}$.
For each $\lambda \in \Lambda_{\chi}$ one defines a baby Verma module $Z_{\chi}(\lambda):= U_{\chi}(\g)\otimes_{U_{\chi}(\mathfrak{b}^{+})} k_{\lambda}$, where $\mathfrak{b}^{+} = \mathfrak{n}^{+} \oplus \mathfrak{h}$ and $k_{\lambda}$ is the one dimensional $U_{\chi}(\mathfrak{b}^{+})$-module defined by $\lambda$. Every irreducible $U_{\chi}(\g)$-module is a factor of a baby Verma module, although in general a baby Verma module can have several irreducible images, and different choices of $\lambda$ can yield the same module \cite[6.7, 6.9]{jan6}.\\
\subsection{The geometric approach to ramification}
\label{mirrum}
Thanks to \cite{mirrum} we get a complete description of the fibres of the map 
$\theta^{*}: \text{Maxspec}(Z)\longrightarrow \text{Maxspec}(Z_0)$ for $\g$ semisimple. In principal 
this yields complete information regarding the ramification index; however we 
would like to have 
a better combinatorial description, for inductive purposes we require the result for reductive $\g$, and, of 
course, we seek representation theoretic consequences.

For $\eta \in \Lambda_{\chi}$ let $W(\eta )$ be the stabiliser of $\eta$ 
and let $W(\eta + \Lambda )$ be the set-wise stabiliser of $\eta +
\Lambda$. Thus $W(\eta)\subseteq W(\eta + \Lambda)$ and 
these are both parabolic subgroups of $W$ by \cite[Proposition 1.15]{humCOX} and 
\cite[Lemma 7]{mirrum}. Consider the partial coinvariant algebra 
\[
C_{\eta} = S(\mathfrak{h}^{(1)})^{W(\eta)} 
\otimes_{S(\mathfrak{h}^{(1)})^{W(\eta+\Lambda)}} k_{\chi_s^p},
\]
where $S(\mathfrak{h}^{(1)})^{W(\eta + \Lambda)} \longrightarrow k_{\chi_s^p}$
is the restriction of $\chi_s^p \in \mathfrak{h}^{(1)*}$. Notice that, since $\eta \in \Lambda_{\chi}$, $W(\eta + \Lambda)$ fixes $\chi_s^p$, so that $C_{\eta}$ is a local ring. It can be shown, \cite[Lemma 9 and (the 
proof of) Theorem 10]{mirrum} that $\dim (C_{\eta}) = [W(\eta +
\Lambda) : W(\eta)]$. Working in the algebra $\Theta(U^G) =
S(\mathfrak{h})^W$, Mirkovi\'{c} and Rumynin prove the first
part of the following result in the case $\g$ is the Lie algebra of a
semisimple algebraic group. In view of Theorem \ref{cent} the
extension is straightforward.
\begin{thm}
Retain as usual the hypotheses and notation of (\ref{liedefs}) and (\ref{az=sing})(A),(B),(C).\\
1. \cite[Theorem 10]{mirrum} Let $\chi\in \g^*$ and let $R_{\chi}$ 
denote a set of representatives of 
$W$-orbits on $W\Lambda_{\chi}$. Then the inverse image of $\chi$ (along 
$\theta^{*}$) is isomorphic via $\Theta$  to the spectrum of $\oplus_{\lambda \in R_{\chi}} 
C_{\lambda}$. That is, 
\begin{equation}
\label{eq:centreman}
Z/\mathfrak{m}_{\chi}Z \quad \cong \quad \bigoplus_{\lambda \in R_{\chi}} 
C_{\lambda}.
\end{equation}
2. For all $\lambda \in \mathfrak{h}^{*}$, $C_{\lambda}$ is a Frobenius algebra. 
\end{thm}
\begin{proof}
2. The key point here is that all $C_{\lambda}$ occur as the centres
of primary components in the Azumaya locus.  This follows from Lemma
\ref{centresred} applied to a semisimple character $\zeta \in \g^{*}$ for which
$\lambda(h)^{p}- \lambda(h^{[p]}) = \zeta(h)^{p}$
holds for all $h\in\h$. For then
$\mathrm{Mat}_{p^N}(C_{\lambda})$ is a direct summand of the
appropriate $U_{\chi}(\g)$, by Corollary \ref{regalgs}. Since
$U_{\chi}(\g)$ is Frobenius by \cite[Proposition 1.2]{frepar1}, it
follows immediately that $C_{\lambda}$ is self-injective. That it is
Frobenius now follows from  \cite[Example, IV.3]{ARS}, since it is commutative and local.
\end{proof}  
\begin{rem}
With our restriction on $p$ the algebras $S(\h)^{W(\eta)}$ can be
realised from integral forms by base change from $\mathbb{Z}$ to
$k$, \cite[9.6]{jan6}. It is therefore straightforward to see that there is an isomorphism  
\[
C_{\eta} \cong S(\h)^{W(\eta)} \otimes_{S(\h)^{W(\eta+\Lambda)}}
k_{\tilde{\chi}_s},
\]
where $\tilde{\chi}_s\in \h^*$ is uniquely defined by
$\tilde{\chi}_s(h_i)=\chi_s(h_i)^p$ for all $1\leq i\leq
r$.
\end{rem}
\subsection{}
\label{cardR}
As is clear from Theorem \ref{mirrum}, the maximal ideals of $Z$ lying over the 
maximal ideal $\mathfrak{m}_{\chi}$ of $Z_0$ are parametrised by $R_{\chi}$. 
This was already known prior to the work in \cite{mirrum} - thus, the baby Verma 
modules $Z_{\chi}(\lambda)$ annihilated by $\mathfrak{m}_{\chi}$ are defined 
for $\lambda$ ranging through the set $\Lambda_{\chi}$ \cite[6.7]{jan6}, 
\cite[Section 10]{hum3}. Conversely, every maximal ideal of $Z$ lying over 
$\mathfrak{m}_{\chi}$ kills such a baby Verma module, and $Z_{\chi}(\lambda)$ 
and $Z_{\chi}(\mu)$ have the same central annihilator if and only if $\lambda$ 
and $\mu$ are in the same $W{\scriptscriptstyle \bullet}$-orbit \cite[Corollary 9.4]{jan6}. Thus in the remainder of this section we shall denote the maximal ideals of $Z$ which lie over a maximal ideal $\mathfrak{m}_{\chi}$ of $Z_0$ by 
\begin{equation*}
\mathfrak{m}_{\chi}Z + \mathfrak{n}_{\eta}Z,
\end{equation*}
where $\eta$ is a representative in $\Lambda_{\chi}$ for a $W{\scriptscriptstyle \bullet}$-orbit in $W\Lambda_{\chi}$ and the maximal ideal $\mathfrak{n}_{\eta}$ of $U^G$ is the inverse image under $(\gamma^{-1}\circ \Theta)^{-1}$ of the maximal ideal of $S(\mathfrak{h})^{W{\scriptscriptstyle \bullet}}$ determined by $\eta$.
\subsection{}
\label{stabs}
We've seen that, thanks to the map $\Theta$ of (\ref{eq:HCclass}), the maximal ideals of $Z$ in 
$\theta^{*^{-1}}(\mathfrak{m}_{\chi})$ are labelled by elements of $\mathfrak{h}^*/W$. 
Theorem \ref{mirrum} shows that $\eta \in \mathfrak{h^*}/W$ yields an unramified point lying over $\mathfrak{m}_{\chi}$ 
if and only if $C_{\eta}= k$; that is, if and only if 
$W(\eta) = W(\eta + \Lambda)$. So we are led to ask: for which $\eta$ does 
$w(\eta) - \eta \in \Lambda$ imply $w(\eta) = \eta$? Since 
$W(\eta + \Lambda)$ is parabolic it is enough to check that $s_{\alpha}( \eta) 
-\eta \in \Lambda $ implies $s_{\alpha} (\eta) =\eta$ for each simple root $\alpha$. 
However it is 
straightforward to check that $s_{\alpha}(\eta) - \eta \in \Lambda$ if 
and only if 
$\eta(h_{\alpha}) \alpha \in \Lambda$, whilst 
$s_{\alpha}(\eta) = \eta$ if and only if $\eta(h_{\alpha})
= 0$.
\begin{thm}
\label{Ugunram}
1. Let $\eta \in \mathfrak{h}^*$. The groups $W(\eta )$ and $W(\eta +
\Lambda)$ are equal if and only if $\eta(h_{\alpha }) \notin \mathbb{F}_p\setminus \{ 0\}$ for all 
simple roots $\alpha$.\\
2. Let $\chi$ be in $\g^{*}$ and let $\eta$ be in  $\Lambda_{\chi}$, with corresponding maximal ideals $\mathfrak{m}_{\chi}$ of $Z_0$ and $\mathfrak{n}_{\eta}$ of $Z_1$. The maximal ideal $\mathfrak{m}_{\chi}Z + \mathfrak{n}_{\eta}Z$ of $Z$ is unramified over $\mathfrak{m}_{\chi}$ if and only if $(\eta + \rho)(h_{\alpha}) \notin \mathbb{F}_p\setminus \{0\}$ for all simple roots $\alpha$.
\end{thm}
\begin{proof}
1. Recall that $Q=\mathbb{Z}\Phi$ is the root lattice and that $Q \subseteq X$ is a free abelian
group. Let $X= X' \oplus X''$ be a decomposition such that
$Q\cap X'' = 0$ and $Q\subseteq X'$ is a
finite extension. By \cite[11.2(1)]{jan6} $Q \cap pX =
pQ$. Thus the map
\[
Q \longrightarrow X \longrightarrow \frac{X}{pX} =
 \Lambda
\]
has kernel $pQ$ and so induces an isomorphism
\begin{equation}
\label{rootisom}
\frac{Q}{pQ} \longrightarrow \frac{X'}{pX'}.
\end{equation}
It therefore follows from (\ref{rootisom}) and the isomorphism $\Lambda\otimes_{\mathbb{F}_p}k \cong
\h^* \cong (Q/pQ \otimes_{\mathbb{F}_p} k)\oplus (X''/pX''\otimes_{\mathbb{F}_p}k)$ that $\eta(h_{\alpha})\alpha \in \Lambda$ if and only if $\eta(h_{\alpha}) \in \mathbb{F}_{p}$, proving the first part of the theorem, in view of the remarks preceding it.\\
\indent 2. Since $\Theta$ involves composition with the winding automorphism $\gamma$, we must apply $\gamma^{-1}$  to Theorem (\ref{mirrum}), so that its conclusions are expressed in terms of representations in $\Lambda_{\chi}$ of $W{\scriptscriptstyle \bullet}$-orbits in $\mathfrak{h}^*$. On doing this, 2. is an immediate consequence of 1.
\end{proof}
\begin{rem}
Some restriction on $p$ is needed for the second
part of the theorem. Indeed if $k$ is a field of characteristic 2 and
$\g = \mathfrak{sl}_2(k)$ then $k[e^2, f^2, h^2 - h]=Z_0 \subset Z =
k[e^2, f^2, h]$. Since the polynomial $h^2 - h - \lambda$ always has two
roots for any $\lambda \in k$ it follows that every point of $\text{Maxspec} (Z)$ is
unramified over $\text{Maxspec} (Z_0)$.
\end{rem}
\subsection{The algebraic approach to ramification}
\label{algram}
By combining Theorem \ref{Ugunram}.2 with Theorems \ref{azunram} and \ref{az=sing}.1 we can deduce a result of Friedlander and Parshall \cite[Theorem 4.2]{frepar2} determining when a baby Verma module $Z_{\chi}(\lambda)$ is projective, as in the equivalence of 1 and 4 below.   For, let $\chi \in \g^*$, (where, as usual, we may assume that $\chi(\mathfrak{n}^+) = 0$), and let $V$ be a simple $U_{\chi}(\g)$-module. Then the annihilator of $V$ in $Z$ is in the Azumaya locus if and only if $V \cong Z_{\chi}(\lambda)$ for some $\lambda \in \Lambda_{\chi}$, since $V$ is an image of such a baby Verma module \cite[Proposition 6.7]{jan6}, and their dimensions are equal precisely when $V$ is annihilated by an Azumaya point of $Z$ \cite[Proposition 3.1]{brogoo}. Thus Theorems \ref{azunram} and \ref{az=sing}.1 imply that an irreducible $Z_{\chi}(\lambda)$ is projective if and only if $\mathfrak{m}_{\chi}Z + \mathfrak{n}_{\lambda}Z$ is unramified over $\mathfrak{m}_{\chi}$, proving the equivalence of 2 and 3 of the following result. The equivalence of 3 and 4 is given by Theorem \ref{Ugunram}.2. Finally, to see that 1 implies 2 follow the argument of \cite[Theorem 4.2]{frepar2} and restrict $Z_{\chi}(\lambda)$ to the local Frobenius subalgebra $U_{0}(\mathfrak{n}^+)$ of $U_{\chi}(\g)$.
\begin{cor}
\label{babyproj}
Let $\chi \in \g^{*}$ with $\chi(\mathfrak{n}^+) = 0$, and let $\lambda \in \Lambda_{\chi}$. Then the following are equivalent.\\
1. $Z_{\chi}(\lambda)$ is projective.\\
2. $Z_{\chi}(\lambda)$ is projective and irreducible.\\
3. $\mathfrak{m}_{\chi}Z + \mathfrak{n}_{\lambda}Z$ is an unramified point of $\mathrm{Maxspec}(Z)$.\\
4. There is no simple root $\alpha$ with $(\lambda +  \rho)(h_{\alpha}) \in \mathbb{F}_{p}-\{0\}$.
\end{cor} 
\subsection{}
\label{morbabyV}
In order to obtain precise combinatorial information we must
reduce to the nilpotent case.
\begin{thm}\cite[Theorem 2]{kacwei1}, \cite[Theorem 3.2]{frepar1}
\label{rednilp}
Let $\chi = \chi_s + \chi_n \in \g^*$ and let $\mathfrak{p} =
\mathfrak{c}_{\g}(\chi_s)\oplus \mathfrak{u}$ be a parabolic
subalgebra of $\g$, having Levi part $\mathfrak{c}_{\g}(\chi_s)$. Then
the functors
\[
F= U_{\chi}(\g)\otimes_{U_{\chi}(\mathfrak{p})}
-:\text{mod }U_{\chi}(\mathfrak{c}_{\mathfrak{g}}(\chi_s))
\longrightarrow \text{mod }U_{\chi}(\g)\]
and
\[
G= (-)^{\mathfrak{u}} : \text{mod }U_{\chi}(\g) \longrightarrow
\text{mod }U_{\chi}(\mathfrak{c}_{\mathfrak{g}}(\chi_s))
\]
are inverse equivalences of categories.
\end{thm}
\begin{rems}
1. We observe that in the above categories baby Verma modules have the
same parameter sets, $\Lambda_{\chi}$. It is straightforward to check
that $F(Z_{\chi}(\lambda)) \cong Z_{\chi}(\lambda)$ for $\lambda \in
\Lambda_{\chi}$ (with the obvious abuse of notation).\\
2. Let $\theta^{-1}(\chi_s) = x_s \in \h$. Recall from Lemma \ref{centresred} that the roots of $\cen{\chi_s}=\cen{x_s}$ are precisely those
for which $\alpha (x_s) = 0$.  By \cite[11.2]{jan6} we have
\[
\theta (h_{\alpha}) = \theta (x_{\alpha})(x_{-\alpha})\alpha,
\]
from which it follows that $\chi_s(h_{\alpha})=(x_s, h_{\alpha}) = 0$ for $\alpha$
appearing in $\cen{x_s}$. Thus we see $\lambda(h_{\alpha})\in
\mathbb{F}_p$ for all roots $\alpha$ in $\cen{\chi_s}$, and $\lambda
\in \Lambda_{\chi}$. \\  
3. Following 2. and (\ref{stabs}) we see immediately that for any
$\lambda\in \Lambda_{\chi}$ the Weyl group of $\cen{\chi_s}$ equals
$W(\lambda + \Lambda)$. In particular $W(\lambda)$ is contained in
the Weyl group of $\cen{\chi_s}$. \\
4. Our hypotheses on $p$ are preserved under passage from $\g$ to $\mathfrak{c}_{\g}(\chi_{s})$, by \cite[6.5]{jan6}. 
\end{rems}
\subsection{}
\label{freparer}
Observe that we have identified two mechanisms whereby a baby Verma module $Z_{\chi}(\lambda)$ can be irreducible: either by virtue of Proposition \ref{regcase} when $\chi$ is regular; or when $\chi$ is arbitrary and (after passing to $\mathfrak{c}_{\g}(\chi_{s})$ by means of Theorem \ref{rednilp}) we have $\lambda = -\rho$ so that $Z_{\chi}(\lambda)$ is a ``module of Steinberg type" and Corollary \ref{babyproj} applies. It is tempting to formulate a result stating, in effect, that these are the only means by which $Z_{\chi}(\lambda)$ can attain irreducibility; and indeed such a result has been stated as \cite[last part of Theorem 4.2]{frepar2}. Unfortunately, this conclusion is false, the problem being that one can construct baby Verma modules of ``mixed type", the simplest such being the following.\\
 
\textbf{Example.} Let $\g = \g_{1} \oplus \g_{2}$ where $\g_{1} \cong \g_{2} \cong \mathfrak{sl}(2)$, and define $\chi \in \g^{*}$ so that $\chi_{1} := \chi|_{\g_{1}}$ is regular nilpotent and $\chi_{2} := \chi|_{\g_{2}} = 0$. Let $\lambda \in \mathfrak{h}^*$, with $\lambda_{1} := \lambda|_{\mathfrak{h}_{1}} \in \Lambda(\chi|_{\mathfrak{h}_{1}})$ and $\lambda_{2}(h_{2}) = \lambda(h_{2}) = -1$, ($h_{2}$ being the Chevalley generator of a Cartan subalgebra of $\g_{2}$). Thus the $U(\g)$-module $Z_{\chi,\lambda} = Z_{\chi_{1},\lambda_{1}} \otimes Z_{\chi_{2},\lambda_{2}}$, being the tensor product of an irreducible $U(\g_{1})$-module with an irreducible $U(\g_{2})$-module, is irreducible; but $\chi = \chi_{n}$ is not regular and $\lambda \neq -\rho$.\\

\indent When $[\g,\g]$ is simple and $\chi|_{[\g,\g]}$ is nilpotent the ``two
mechanisms" result is true, and the proof given in \cite{frepar2}
works for this case; so using Theorem \ref{rednilp} the interested
reader can easily formulate the correct result in general. Let us remark also that one can approach this result geometrically
using Theorem \ref{az=sing}.1. In particular one can describe the
simple $U_{\chi}(\g)$-modules of maximal dimension by following the proof of
\cite[Theorem 2]{kry} using the description of the centre obtained
earlier.  We shall
prove one direction of the analogous result for quantum groups in
Theorem \ref{simpqg}.\\

\subsection{}
Theorem \ref{rednilp} allows us to describe all the unramified maximal ideals of $Z$
lying over $\mathfrak{m}_{\chi}$ for any $\chi\in \g^*$.
\begin{prop}
Let $\chi\in \g^*$ be such that $\chi(\n^+)=0$, so that $\h \subseteq
\mathfrak{c}_{\g}(\chi_{s})$. Let $\pi : \h^* \rightarrow
([\cen{\chi_s},\, \cen{\chi_s}]\cap \h )^*$ be projection. Let $s=dim (ker
\pi)$. Then there are exactly
$p^s$ unramified maximal ideals of $Z$ lying over $\mathfrak{m}_{\chi}$. They
are given by $\mathfrak{m}_{\chi}Z+\mathfrak{n}_{\lambda}Z$ where $\lambda \in
\pi^{-1}(\pi(-\rho))$.
\end{prop}
\begin{proof}
Suppose first that $\cen{\chi_s} =\g$. By Corollary \ref{babyproj},
$\mathfrak{m}_{\chi}Z+\n_{\lambda}Z$ is unramified if and only if
$(\lambda+\rho)(h_{\alpha}) \notin \mathbb{F}_p\setminus \{
0\}$ for all simple roots $\alpha$, so the second remark following Theorem \ref{rednilp} forces $(\lambda +\rho)(\h_{\alpha})=0$ if $\mathfrak{m}_{\chi}Z+\n_{\lambda}Z$ is unramified. So we have that $\pi (\lambda) = \pi(-\rho)$
showing that the unramified maximal ideals of $Z$ lying above
$\mathfrak{m}_{\chi}$ are of the required form. Now suppose $\lambda$ and $\nu$
are in $\pi^{-1}(\pi(-\rho))$ and give rise to the same maximal ideal of $Z$. In other words
$Z_{\chi}(\lambda)\cong Z_{\chi}(\nu)$ and so $\lambda$ and $\nu$ are
in the same $W{\scriptscriptstyle \bullet}$-orbit. But the stabiliser of $\lambda$ is all
of $W$, which implies $\lambda =
\nu$. Thus the proposition is true in this case.

To prove the general case apply Theorem \ref{rednilp} together with the remarks
following it.
\end{proof}
\subsection{The Azumaya locus}
\label{azlocsec}
The proof of \cite[Theorem 4.10]{brogoo} relies on the fact that when
$\chi$ is regular all the
simple $U_{\chi}(\g)$-modules have dimension $p^N$, the PI degree of
$U(\g)$. We will prove the converse in this paragraph: if all simple
$U_{\chi}(\g)$-modules have dimension $p^N$ then $\chi$ is regular. We
will also give in the next paragraph a combinatorial description of the algebra
$U_{\chi}(\g)$ in this case.
\begin{prop}
Let $\chi\in\g^*$. Then $\mathfrak{m}_{\chi}$ is fully Azumaya if and only if $\chi$ is
regular.
\end{prop}
\begin{proof}
Suppose first that $\chi$ is regular. Then, by \cite[Theorem 3.10]{pre2} (noting that $p \neq 2$ by hypothesis),
$p^{\frac{1}{2}\dim(G.\chi)}=p^N$ divides the dimension of all the
simple $U_{\chi}(\g)$-modules, as required.

Conversely, suppose that $\chi$ is not regular. Since the simple
modules of maximal dimension are necessarily baby Verma
modules we see from the remarks following Lemma \ref{centresred} and Theorem \ref{rednilp} that we may assume that
$\cen{\chi}=\g$. Therefore $\chi_n$ is not regular nilpotent. There
must exist a simple root of $\g$, say $\alpha$, such that
$\chi_n(x_{-\alpha}) = 0$ (for otherwise $\g$ is a torus, whence
$\chi$ is regular). There exists $\lambda \in \Lambda$ such that
$(\lambda+\rho)(h_{\alpha})\neq 0$. Now the argument of \cite[Theorem
4.2]{frepar2} is valid in this context and shows that
$Z_{\chi}(\lambda)$ is not simple. Indeed let $\mathfrak{p}$ be the
minimal parabolic subalgebra of $\g$ with Levi subalgebra spanned by
$x_{\alpha}, x_{-\alpha}$ and the elements of $\h$. Then $\chi$
restricted to $\mathfrak{p}$ is zero on $x_{\alpha}$ and $x_{-\alpha}$
so $\mathfrak{sl}_2$ representation theory tells us that
$U_{\chi}(\mathfrak{p})$ has a simple module, $M$, of dimension less than $p$
whose highest weight is $\lambda$, \cite[Lemma 7.2]{HUMLA}. Thus there is a non-zero
homomorphism from $Z_{\chi}(\lambda)$ to
$U_{\chi}(\g)\otimes_{U_{\chi}(\mathfrak{p})} M$. Since the right hand
side has dimension equal to $p^{N-1}\dim M < p^N$ we are done.
\end{proof}
\subsection{}
\label{regcase}
We will now describe the algebras $U_{\chi}(\g)$ for $\chi \in \g^*$
regular. Without loss of generality $\chi(\n^+) = 0$. Recall from
Lemma \ref{centresred} that $\mathfrak{c}_{\g}(\chi_s)$ is the
subalgebra of $\g$ spanned by $\h$ and the vectors $\{x_{\alpha}:
\chi_s(h_{\alpha}) =0\}$, a
Levi subalgebra of $\g$. Let $\lambda \in \Lambda_{\chi}$. Then note
that, by the proof of Theorem \ref{Ugunram} and Remark \ref{morbabyV}.2, $W(\lambda + \Lambda)$ is $W'$, the
Weyl group of $\cen{\chi_s}$, and so $W(\lambda) = W'(\lambda)$, the
stabiliser of $\lambda$ in $W'$.
 
\begin{prop}
Keep the above notation. In particular $\chi\in \g^*$ is regular with $\chi(\mathfrak{n}^+) = 0$ and
$W'$ is the Weyl group of $\cen{\chi_s}$. Let $\pi : \h^* \longrightarrow
(\h\cap [\cen{\chi_s},\, \cen{\chi_s}])^*$ be projection.  Then
\[
U_{\chi}(\g) \cong \mathrm{Mat}_{p^N}\left(
\bigoplus_{\lambda\in \Lambda_{\chi}}
S(\h^{(1)})^{W'(\lambda)} \otimes_{S(\h^{(1)})^{W'}}k_{\chi_s}\right) .
\]
Moreover if $\lambda, \eta\in \Lambda_{\chi}$ are such that
$\pi(\lambda) = \pi(\eta)$ then 
\[
S(\h^{(1)})^{W'(\lambda)}
\otimes_{S(\h^{(1)})^{W'}}k_{\chi_s}\cong S(\h^{(1)})^{W'(\eta)}
\otimes_{S(\h^{(1)})^{W'}}k_{\chi_s}.
\]
\end{prop}
\begin{proof}
As we have seen $W(\lambda + \Lambda) = W'$ and $W(\lambda ) =
W'(\lambda)$. Therefore
\[
C_{\lambda} = S(\h^{(1)})^{W'(\lambda)}
\otimes_{S(\h^{(1)})^{W'}}k_{\chi_s}.
\]
The first isomorphism therefore follows from Corollary \ref{regalgs}, Theorem
\ref{mirrum} and Proposition \ref{azlocsec}. The assumption that $\pi(\lambda) = \pi
(\eta)$ implies that $\lambda(h_{\alpha}) = \eta
(h_{\alpha})$ for all $\alpha$ a root of $\cen{\chi_s}$ and so
$W(\lambda) = W(\eta)$ as required.
\end{proof}
\subsection{}
\label{lieblocks}
\label{Muller2}
We would like to be able to describe the centre $Z(U_{\chi})$ of $U_{\chi}$, and in particular to compare it with the algebra $Z_{\chi} = Z/\mathfrak{m}_{\chi}Z$ described by the result of Mirkovi\'{c} and Rumynin, Theorem \ref{mirrum}. Thus, we can ask:\\

A. Is the natural homomorphism $\Psi$ from $Z_{\chi}$ to $Z(U_{\chi})$ injective? Equivalently, is $\mathfrak{m}_{\chi}U \cap Z = \mathfrak{m}_{\chi}Z$?\\

B. Is $\Psi$ onto?\\

 For regular $\chi$ the answer to Question A is ``yes", \cite[Lemma 11]{mirrum}; this result can also be read off from Corollary \ref{regalgs}, and this conclusion has recently been extended to arbitrary $\chi$ in $\g^{*}$ by Premet \cite{precom}. The answer to question B is ``no", as is demonstrated by the following example which was brought to our attention by Sasha Premet. We are grateful to him for allowing us to include it here.\\
\begin{example}[Premet]
 Continue with the hypotheses and notation of (\ref{liedefs}), but
 assume $G$ is semisimple. Then there is a (nilpotent) element of the centre of the restricted enveloping algebra $\overline{U}$ of $\g$ which is not in the image of the centre $Z(\g)$ of $U(\g)$ under the canonical homomorphism from $U(\g)$ onto $\overline{U}$.
\end{example}
\begin{proof} 
Let $\Theta : U(\g)^G \longrightarrow U(\h)^W$ be the isomorphism of
Theorem \ref{centreU}. Giving $U(\g)$ its natural filtration we can see, by passing to
integral forms as in \cite[9.6]{jan6}, that $\Theta$ is a filtration
preserving map, where $U(\h)^W$ has its natural graded structure. As
we have seen in (\ref{adjquot}) the algebra $U(\h)^W$ has algebraically
independent homogeneous generators $T_1,\ldots ,T_r$ of degree
$d_1,\ldots ,d_r$ respectively. Hence, if we let
$Y_i=\Theta^{-1}(T_i)$ for all $i$, we see that $U(\g)^G$ also has
algebraically independent generators of degrees $d_1,\ldots ,d_r$.
By \cite[Theorem 3.1]{veld}, which we have already noted is valid
thanks to the earlier results of this chapter,
$Z=Z(g)$ is freely generated as a $Z_0$-module by $\{Y_1^{i_{1}}
\ldots Y_r^{i_{r}} : 0 \leq i_j \leq p-1, 1 \leq j \leq r \}$. Set
$\overline{Z}$ to be the image of the centre $Z$ of $U$ in
$\overline{U}$, so that $\overline{Z} = k_0\otimes_{Z_1\cap
Z_0}Z_1$. Then we deduce from the foregoing that the maximum degree of
an element of $\overline{Z}$ is $\sum_{i=1}^r(p-1)d_i$. As before we have $\sum_{i=1}^rd_i = N+r$, by
\cite[Theorem 3.9]{humCOX}.

We now exhibit an element of $Z(\overline{U})$ with degree $(p-1)(2N + r)$,
thus proving that the map from $\overline{Z}$ to $Z(\overline{U})$ is not
surjective as long as $N > 0$. Recall that $\g$ acts on $U$ and $S(\g)$ by the adjoint
action, and this induces an action of $\g$ on $\bar{U}$ and on
$\overline{S(\g)}$, where the latter notation has the obvious
meaning. By \cite[Theorem 2.1]{frepar7} there is a filtration preserving isomorphism of $\g$-modules
\begin{equation*}
\beta : \overline{U} \longrightarrow \overline{S(\g)}
\end{equation*}
induced by the Mil'ner map. Up to scalars
there is a unique element $z$ of degree $(p-1)(2N + r)$ in
$\overline{S(\g)}$, and the space it spans must therefore be
$\g$-stable, and, since $\g$ is semisimple, in fact trivial. The $\g$-equivariance of $\beta$ shows that $\beta^{-1}(z)$ is the central element of $\overline{U}$ of large degree which we wanted to find.
\end{proof}
\begin{rem}
The above proof carries over to the case $G$ reductive and $p$ very
good (that is $p$ is good and $(p,n+1)=1$ if $\mathcal{D}G$ has a
component of type of $A_n$). To see this we need only observe that
$\g$ decomposes as $\text{Lie}(\mathcal{D}G)\oplus \mathfrak{z}$,
where $\mathfrak{z}$ is the centre of $\g$.
\end{rem}
\subsection{} 
\label{gblocks}
In fact as shown in (\ref{Blocks}) it is with the nilpotent elements of $Z(U_{\chi})$ that the failure of $\Psi$ to be surjective lies. Applying this result in the present enveloping algebras setting yields a description of the blocks of $U(\g)$ which gives confirmation of a conjecture of Humphreys \cite[Section 18]{hum3}, who noted there the case of standard Levi type, generalising his 1971 result for the case $\chi=0$ \cite{hum20}. The following is an immediate consequence of (\ref{cardR}) and Corollary \ref{Muller}.
\begin{thm}
Let $\chi \in \g^{*}$.  The blocks of $U_{\chi}(\g)$ are in natural bijection 
with the $W{\scriptscriptstyle \bullet}$-orbits in $W\Lambda_{\chi}$.
\end{thm} 
\indent In the light of this result we shall label the blocks of $U_{\chi}(\g)$ by $W{\scriptscriptstyle \bullet}$-orbit representatives from $W{\scriptscriptstyle \bullet}\Lambda_{\chi}$. Thus for a $W{\scriptscriptstyle \bullet}$-orbit representative $\lambda$ from $W{\scriptscriptstyle \bullet}\Lambda_{\chi}$ we write $\mathcal{B}_{\chi,\lambda}$ for the block of $U_{\chi}(\g)$ whose simple
modules are annihilated by the maximal ideal $\mathfrak{m}_{\chi} + \mathfrak{n}_{\lambda - \rho}$ of $Z$. One can immediately read off from Theorems \ref{azunram} and \ref{mirrum}, and Corollary \ref{algram} the equivalence of the following statements:
\begin{enumerate}
\item $\mathcal{B}_{\chi,\lambda}$ is simple Artinian;
\item $\mathfrak{m}_{\chi} + \mathfrak{n}_{\lambda - \rho}$ is unramified over $\mathfrak{m}$;
\item no simple root $\alpha$ has $\lambda(h_{\alpha}) \in \mathbb{F}_{p} - \{0\}$;
\item $W(\lambda) = W(\lambda + \Lambda)$.
\end{enumerate}
\indent Naturally one next considers blocks of finite representation type. These have been determined when $\chi$ is standard Levi and $\g$ is classical simple satisfying hypotheses (A), (B) and (C) by \cite{nakpol}; we show in the next section how to obtain their result, extended to reductive $\g$ so that inductive methods based on Theorem \ref{morbabyV} can be used to deal with non-nilpotent characters, using ramification methods.  
\subsection{Blocks of finite representation type}
\label{reptype}
Let $\Delta'$ be a subset of $\Delta$ and let $W' = < s_{\alpha} : \alpha \in \Delta '>$ be a parabolic subgroup
of $W$. Recall that there is a
distinguished set of coset representatives for $W'$ in $W$, called the minimal coset
representatives. These are the set $\{ w\in W:
\ell(ws_{\alpha})>\ell(w) \text{ for all }\alpha\in \Delta'\}$, where
$\ell(w)$ denotes the length of $w\in W$, \cite[Proposition
1.10]{humCOX}. In case $W' = W(\lambda)$ for some $\lambda \in \h^*$
we will write $W^{\lambda}$ for this set of representatives. 
\begin{lem}
Let $\chi= \chi_n$ be nilpotent and let $\lambda \in \Lambda$. Then
$C_{\lambda}$ is a graded algebra with Poincar\'{e} series
\[
P(C_{\lambda}, t) =  \sum_{w\in W^{\lambda}} t^{\ell (w)}.
\]
In particular, $C_{\lambda}$ is uniserial only if no two distinct members of $W^{\lambda}$ have the same length.
\end{lem}
\begin{proof}
Since $\chi_s = 0$ it is clear from the definition in (\ref{mirrum}) that $C_{\lambda}$ is graded. By \cite[Section 3]{dem}
for every $w\in W$ there is an operator, $D_w$, on
$S(\h)$ which lowers degree by $\ell(w)$.  By
\cite[9.6]{jan6} the characteristic of the field, $p$, is not a
torsion prime (in the sense of \cite{dem}) and by hypothesis it is not 2. Thus there exists an
element $a\in S^N(\h)$ such that the elements $D_w(a)$ give a basis
for $S(\h)$ over $S(\h)^W$, \cite[Th\'{e}or\`{e}me 2 and Corollaire]{dem}. 

Let $w\in W^{\lambda}$. Then we claim that $D_{w^{-1}w_0}(a)$ is
$W(\lambda)$-invariant. To see this we use the fact that $D_wD_{w'} =
0$ unless $\ell(ww') = \ell(w) + \ell(w')$. By the definition of
$W^{\lambda}$ we see then that $D_{s_{\alpha}}D_{w^{-1}w_0}=0$ for all
simple reflections $s_{\alpha}$ in $W(\lambda)$. Thus, by \cite[Section 3]{dem},
we have for all simple reflections $s_{\alpha} \in W(\lambda)$
\[
0 = (\alpha\otimes 1)D_{s_{\alpha}}(D_{w^{-1}w_{0}}(a)) = D_{w^{-1}w_{0}}(a) -s_{\alpha}D_{w^{-1}w_{0}}(a),
\]
as claimed.

From the previous paragraphs we see that the elements $D_{w^{-1}w_0}(a)$,
for $w\in W^{\lambda}$, span a subspace of
$S(\h)^{W(\lambda)}\otimes_{S(\h)^W}k_0$ of dimension
$[W:W(\lambda)]$. Moreover, by Remark \ref{mirrum}, $C_{\lambda} \cong
S(\h)^{W(\lambda)}\otimes_{S(\h)^W}k_0$ and, as we have noted already, it follows from
\cite[Lemma 9 and Theorem 10]{mirrum} that $C_{\lambda}$ has dimension
$[W:W(\lambda)]$. Thus the elements $D_{w^{-1}w_0}(a)$ with $w\in
W^{\lambda}$ yield a basis for
$C_{\lambda}$ and since $\text{deg}(D_{w^{-1}w_0}(a)) = N-
\ell(w^{-1}w_0) = \ell(w)$, the Poincar\'{e} series is as claimed. For
the last part, simply notice that if a finite dimensional
$\mathbb{N}$-graded algebra has a homogeneous component of dimension
greater than 1 then it cannot be uniserial.
\end{proof}
Recall a nilpotent character of $\g$ is called \textit{standard Levi} if it is
regular nilpotent in some Levi subalgebra of $\g$. The following result generalises \cite[Theorem 4.2]{nakpol} which
determines when a block of a reduced enveloping algebra corresponding to
a nilpotent character of standard Levi type is of finite
representation type in the case when $\g$ is the Lie algebra of simple,
simply-connected algebraic group and the characteristic of $k$ is very
good.  
\begin{cor}
Keep hypotheses (\ref{liedefs}) and (\ref{az=sing})(A), (B) and
(C). Let $\chi= \chi_s+\chi_n\in \g^*$ and let $\lambda \in \Lambda_{\chi}$.
\newline
1. The block $\mathcal{B}_{\chi,\lambda}$ has finite representation
type and is not semisimple only if the rank of $W(\lambda)$ is one less than the rank of $W(\lambda+\Lambda)$
and on the connected Coxeter graph on which these groups differ one of the following holds:
\newline
\indent
(i) $W(\lambda + \Lambda)$ is of type $A_n$ and $W(\lambda)$ is of type
$A_{n-1}$ ;
\newline
\indent
(ii) $W(\lambda + \Lambda)$ is of type $B_n$ and $W(\lambda)$ is of type
$B_{n-1}$ ;
\newline
\indent
(iii) $W(\lambda + \Lambda)$ is of type $G_2$ and $W(\lambda)$ is of
type $A_1$.
\newline
\noindent
(Here we take $A_0 = \emptyset$ and $B_1 = A_1$.)
\newline
2. Assume $\chi_n$ is of standard Levi type and the block
$\mathcal{B}_{\chi,\lambda}$ has finite representation type. Then
$Z_{\chi}(\lambda)$ is the unique simple $\mathcal{B}_{\chi,\lambda}$-module.
\newline
3. If $W(\lambda)$ has rank one less than $W(\lambda+\Lambda)$, one of
(i), (ii) and (iii) holds and $Z_{\chi}(\lambda)$ is the unique simple $\mathcal{B}_{\chi,\lambda}$-module then
$\mathcal{B}_{\chi,\lambda}$ has finite representation type.
\end{cor}
\begin{proof}
First note that by Remark \ref{morbabyV}.3 the groups $W(\lambda)$ and
$W(\lambda+\Lambda)$ lie inside the Weyl group of $\cen{\chi_s}$. It
follows from Theorem \ref{morbabyV} that we can (and will) assume without loss of
generality that $\g = \cen{\chi_s}$. 

Suppose that $A$ is a finite dimensional uniserial algebra (with a
unique simple module). Then $A$ is Morita equivalent to a truncated
polynomial ring $k[X]/(X^n)$ for some $n\in \mathbb{N}$. To see this
let $P$ be the projective cover of the simple $A$-module and consider
$\text{End}_A(P)$. As $P$ is uniserial with all composition factors
isomorphic so it is easy to check that 
\[
\text{End}_A(P) = \frac{k[X]}{(X^n)},
\]
where $n$ is the Loewy length of $P$.

Suppose that $\mathcal{B}_{\chi,\lambda}$ has finite representation
type. By \cite[Theorem 3.2]{far1} $\mathcal{B}_{\chi,\lambda}$ is
uniserial and has a unique simple module. The above paragraph shows
that the centre of $\mathcal{B}_{\chi,\lambda}$ is a truncated
polynomial ring and so, by Premet's positive answer \cite{precom} to \ref{lieblocks}A, $C_{\lambda}$ is a subalgebra of a truncated
polynomial ring $T$. Moreover, $C_{\lambda}$ is a Frobenius algebra, by Theorem \ref{mirrum}.2. Let $I$ be an ideal of $C = C_{\lambda}$. By \cite[Proposition XIV.2.2(ii)]{sten}, $I = \mathrm{Ann}_{C}(J)$ for some ideal $J$ of $C$, and hence $I = \mathrm{Ann}_{T}(JT) \cap C$. In particular $C_{\lambda}$ is a uniserial algebra (and in fact  a truncated polynomial algebra). Following Lemma \ref{reptype} we see that $C_{\lambda}$ is
 uniserial
 only if the set of minimal coset representatives $W^{\lambda}$
 has at most one element of any given length. In particular for
 $C_{\lambda}$ to be uniserial it follows
 that the rank of $W(\lambda)$ must be no less than the rank of $W$
 minus one. On the other hand if $W(\lambda) = W$ then $C_{\lambda}$
 is semisimple, so we need only consider the case when the rank of
 $W(\lambda)$ is one less than the rank of $W$. 

Let $\Gamma$ be the Coxeter graph associated with $W$ and suppose that $s_{\alpha_i}$ is the simple reflection not in
$W(\lambda)$. There are a few cases to consider. Firstly suppose that there is a subdiagram of $\Gamma$ of the
following form, where we label the nodes from left to right as $i-1$,
$i$ and $i+1$,
\[
\xymatrix{
\cdots {\circ} \ar@{-}[r] & {\circ} \ar@{-}[r] & {\circ}
 \cdots
}
\]
Then both $s_{i-1}s_i$ and $s_{i+1}s_i$
are elements of length two in $W^{\lambda}$, implying that $C_{\lambda}$ is not
uniserial. Secondly, if $i$ occurs as one of the endpoints in the subdiagram of
$\Gamma$ below, then label the nodes from top to bottom and
left to right by $i, i+1, i+2$ and $i+3$,
\[
\xymatrix{
 & \circ \ar@{-}[d] &  \\
\cdots \circ \ar@{-}[r] & \circ \ar@{-}[r] &\circ \cdots}
\]
Then both
$s_{i+1}s_{i+2}s_i$ and $s_{i+3}s_{i+2}s_i$ are elements of length
three in 
$W^{\lambda}$, implying that $C_{\lambda}$ is not uniserial (this
argument applies to type $D$ and type $E$). Finally
assume that $i$ occurs in the subdiagram of $\Gamma$ below where we
label the nodes from left to right by $i$, $i+1$ and $i+2$,
\[
\xymatrix{
\circ \ar@2{-}[r] & \circ \ar@{-}[r] & \circ \cdots }
\]
Then both $s_{i+2}s_{i+1}s_{i}$ and $s_is_{i+1}s_i$ are
elements of length three in $W^{\lambda}$, showing $C_{\lambda}$ is not uniserial (the
same argument works for $F_4$). This proves the first part of the corollary.

Since we are assuming $\g = \cen{\chi_s}$ we see by Remark \ref{morbabyV}.2 that $\chi_s(h_{\alpha})=0$ for
all roots $\alpha$. Thus there exists $\mu \in \Lambda_{\chi}$ such
that $\mu(h_{\alpha}) = 0$ for all roots $\alpha$. There is a
bijection
\[
\tau_{\mu} : \Lambda_{\chi} \longrightarrow \Lambda_{\chi_n}=\Lambda,
\]
defined by $\tau_{\mu}(\lambda) = \lambda -\mu$ for $\lambda\in\Lambda_{\chi}$. Moreover
$\tau_{\mu}$ is $W$-equivariant since $s_{\alpha}(\mu) = \mu$ for all
roots $alpha$. We can combine \ref{freparer} (in particular the refined
version of \cite[Theorem 4.2]{frepar2}) together with \cite[Proposition
4.1]{nakpol} (which applies since we can replace $\lambda$ by
$\tau_{\mu}(\lambda)$ as above and therefore assume that $\chi$ is
nilpotent) to see that if $\chi_n$ is of standard Levi type and
$Z_{\chi}(\lambda)$ is not simple then $\mathcal{B}_{\chi , \lambda}$
has more than one simple module. Thus, by \cite[Theorem 3.2]{far1},
$\mathcal{B}_{\chi,\lambda}$ is not of finite representation
type. This proves the second claim. 

Finally, assume the conditions of 3. hold. Then by Lemma \ref{regcase} we have
\[
\mathcal{B}_{\chi,\lambda} \cong \text{Mat}_{p^N}(C_{\lambda}).
\]
In particular we see that the representation type of
$\mathcal{B}_{\chi,\lambda}$ depends only on $C_{\lambda}$. Therefore
we can replace $\chi$ by a regular nilpotent element of $\g^*$, say $\chi'$, and
$\lambda$ by $\tau_{\mu}(\lambda)$ (as above) and deduce that
$\mathcal{B}_{\chi,\lambda}$ has finite representation type if and
only if $\mathcal{B}_{\chi',\tau_{\mu}(\lambda)}$ has finite
representation type. Now the (short) argument
in \cite[Theorem 4.2]{nakpol} shows indeed that $\mathcal{B}_{\chi',
\tau_{\mu}(\lambda)}$ has finite representation type. This proves the final claim and so completes the corollary.
\end{proof} 
\section{Quantised enveloping algebras at roots of unity}
\subsection{}
\label{notqg}
Let $\ep\in \C$ be a primitive root of unity of order $\ell$ and assume 
throughout that $\ell$ is odd and prime to $3$ if $\g$ has a component of type 
$G_2$. Throughout this section $k$ is the complex field $\mathbb{C}$. We continue to use whenever relevant the notation for weights, Weyl group and so on introduced in  (\ref{liedefs}). Let $\UE$ be the simply connected quantised enveloping algebra of $\g$ at 
the root of unity $\ep$. In particular, this means that we should
consider $\g$ as the Lie algebra of the \textit{simply connected}
group $G$. Let $Z_0=Z_0^-\otimes Z_0^0\otimes Z_0^+$ be the 
central subalgebra of $\UE$ generated by the $\ell^{\text{th}}$-powers of the canonical generators and let 
$Z_1$ be the central subalgebra of $\UE$ obtained from the centre of the generic 
quantised enveloping algebra by specialisation. For details of these algebras see \cite[Section 21]{decpro3}.

Let $\UO$ be the subalgebra of $\UE$ generated by the elements $K_{\lambda}$ for 
$\lambda \in X$. This is just the ring of functions of the maximal torus, $T$, 
of $G$, and $\UO$ is the complex group algebra of the free Abelian
group $\langle K_{\varpi_{i}} : 1 \leq i \leq r \rangle$, where
$\{\varpi_{i} : 1 \leq i \leq r \}$ are a set of fundamental weights
of $\g$, orthogonal to the simple roots $\{\alpha_{i} : 1 \leq i \leq
r \}$. Let $\Gamma$ be the group of homomorphisms from $X$ to
$\mathbb{Z}_2=\{\pm 1\}$, so $\Gamma$ is an elementary abelian two group of rank $r$. There is an action of $\Gamma$ on $U_0$ 
given by $\mu\cdot K_{\lambda} = \mu(\lambda)K_{\lambda}$, for $\mu \in \Gamma$. The fixed 
ring, $(\UO)^{\Gamma}$, is again the ring of functions of a maximal torus of 
$G$, which we can identify with $T$. Under this identification the inclusion 
$(\UO)^{\Gamma}\subset \UO$ induces the morphism $\sigma : T \longrightarrow T$, 
where $\sigma (t) = t^2$. Similarly $Z_0^0$ can be considered as the ring of 
functions of $T$ and, in this case, the inclusion $Z_0^0\subset \UO$ induces the 
morphism $F: T\longrightarrow T$, where $F(t) = t^{\ell}$. We shall sometimes write $U_{0,\ell}$ for the subalgebra $Z^{0}_{0}$ of $Z_0$.

Recall the Harish-Chandra map
\begin{equation}
\label{eq:HCmap}
\begin{CD}
\psi : Z_1 @> \pi_{\UO} >> \UO @> \gamma >> \UO,
\end{CD}
\end{equation}
where $\pi_{\UO}$ is projection onto $\UO$ obtained using the triangular decomposition $\UE = U^+ \otimes U_0 \otimes U^-$ \cite[Theorem 19.1]{decpro3} and $\gamma (K_{\lambda}) = 
\ep^{-(\rho, \lambda)}K_{\lambda}$.  By \cite[Section 21]{decpro3} $\psi$ is an 
injective map whose image is $\UGW$.

\begin{thm}\cite[Sections 19,20,21]{decpro3}
\newline
\noindent
1. There is an identification of $\mathrm{Maxspec}(Z_0^{\pm})$ with $U^{\pm}$. 
Let $G^0=B^-B^+$ denote the big cell of $G$. The morphism
\[
\pi: \mathrm{Maxspec}(Z_0) \cong (U^-  \times U^+)\rtimes T \longrightarrow G^0 \subset 
G,
\]
taking $(u_- u_+,t)$ to $u_-^{-1}t^2u_+$ is an unramified covering of $G^0$ of 
degree $2^r$.
\newline
2. Each symplectic leaf of $\mathrm{Maxspec}(Z_0)$ contains an element
$\chi =\chi_u\chi_s \in  
U^-\rtimes T=B^-$ such that $\pi(\chi_u)\pi(\chi_s)$ is the Jordan
decomposition of $\pi(\chi)$ in $G$.
\newline
3. There is an algebra isomorphism
\[
Z(\UE ) \cong Z_0 \otimes_{Z_0\cap Z_1} Z_1.
\]
4. The image of  $Z_0\cap Z_1$ under the monomorphism $\psi$ is $\ULGW$.\\
5. As a $(Z_0  \cap Z_1)$-module, $Z_1$ is free of rank $\ell^{r}$, and $Z(\UE)$ is a free $Z_0$-module of rank $\ell^r$.
\end{thm}

\subsection{}
\label{hearmysong}
Let $t\in T$ and let $C_{G}(t)$ be the centraliser of $t$ in
$G$. Then, by \cite[Theorems 2.2 and 2.11]{humCC}, $C_G(t)$ is a
connected, reductive group which is generated by $T$ together with
those unipotent subgroups $U_{\beta}$ for which $\beta (t)=1$. In
particular, if $W(t)$ denotes the Weyl group of $C_G(t)$ then $W(t) =
<s_{\beta} : \beta (t)=1>$. This is not necessarily a parabolic
subgroup of $W$. By \cite[Theorems 2.12 and 2.15]{humCC}, however,
the root system $\Phi_t = \{ \beta : \beta(t) = 1\}$ of $C_G(t)$ has a
basis which is $W$-conjugate to a proper subset of $\tilde{\Delta}=
\Delta\cup \{-\alpha_0\}$, where $\alpha_0$ is the highest root of
$\Phi^+$. 
\begin{lem}\cite[Lemma 6.1(d)]{deckac100}
Suppose $\chi=\chi_u\chi_s$ such that $\pi(\chi_u)\pi(\chi_s)$ is the
Jordan decomposition of $\pi(\chi)$. Then if
$\chi_s(K_{\beta}^{2\ell})\neq 1$ for some $\alpha\in \Phi^+$ then
$\chi_u(F_{\beta}^{\ell}) = 0$.
\end{lem}
\begin{proof}
For $\beta\in \Phi^+$ let $y_{\beta} = (\ep^{(\beta ,\beta)/2} -
\ep^{-(\beta ,\beta )/2})^{\ell}F_{\beta}^{\ell}$. Then by definition,
\cite[Section 4]{deckacpro4},
\[
\pi(\chi_u) = \text{exp}(\chi_u(y_{\beta_N})f_{\beta_N})\ldots
\text{exp}(\chi_u(y_{\beta_1})f_{\beta_1}).
\]
On the other hand since $\pi(\chi_u)\in C_G(\pi(\chi_s))$ and the
decomposition above is unique (for that ordering of positive roots) we
must have $\chi_u(y_{\beta})=0$ for all $\beta \in \Phi^+$ such that
$\beta (\pi(\chi_s)) \neq 1$. By definition $\beta (\pi(\chi_s)) =
\chi_s (K_{\beta}^{\ell})^2$, proving the lemma.
\end{proof}
\subsection{} 
\label{size}
We now define and study the quantum analogue of the algebras $C_{\eta}$ of (\ref{mirrum}).
For $t\in T$ define the algebra
\[
D_{t} :=  U_{0,\ell} ^{\Gamma . W(t)} \otimes_{U_{0,\ell}^{\Gamma 
.W(t^{\ell})}} \C_{t^{\ell}},
\]
where $U_{0,\ell}^{\Gamma . W(t^{\ell})}\longrightarrow \C_{t^{\ell}}$
is the restriction of the map $U_{0,\ell}^{\Gamma} \longrightarrow
\C_{t^{\ell}}$ which sends $K_{\lambda}^{2\ell}$ to $\lambda (t)^{\ell}$.
\begin{lem}
Let $t\in T$. Then there is an isomorphism
\[
D_t \cong U_0^{W(t)} \otimes_{U_0^{W(t^{\ell})}} \C_{t^{\ell}},
\]
where $U_0^{W(t^{\ell})}\longrightarrow \C_{t^{\ell}}$ is evaluation at 
$t^{\ell}$. Thus $D_t$ is a local algebra. Moreover for any $w\in W$ there is an isomorphism $D_t \cong 
D_{w(t)}$.
\end{lem}
\begin{proof}
The algebras $U_0$ and $U_0^{\Gamma}$ (respectively $U_0$ and $U_{0,\ell}^{\Gamma}$) are 
isomorphic via the map sending $K_{\lambda}$ to $K_{\lambda}^2$ (respectively 
$K_{\lambda}$ to $K_{\lambda}^{2\ell}$). It is clear that these isomorphisms are 
$W$-equivariant. Moreover $U_{0,\ell}^{\Gamma.W(t^{\ell})} \longrightarrow \mathbb{C}_{t^\ell}$ is the restriction 
of the map on $U_{0,\ell}^{\Gamma}$ which sends $K_{\lambda}^{2\ell}$ to $K_{\lambda}^{2\ell}(t^{\ell}) = 
\lambda (t)^{\ell}$. Under the above isomorphism this becomes the map sending 
$K_{\lambda}$ to $\lambda(t)^{\ell}$ so it follows that
\[
D_t \cong U_0^{W(t)}\otimes_{U_0^{W(t^{\ell})}} \C_{t^{\ell}}.
\]

The second statement is now clear. For the final statement let $s=w(t)$ and note that $W(s) = wW(t)w^{-1}$. Thus 
there is an isomorphism from $U_0^{W(t)}$ to $U_0^{W(s)}$ taking $x$ to $w.x$. 
Under this isomorphism the map $U_0^{W(t^{\ell})}\longrightarrow \C_{t^{\ell}}$ 
obtained by evaluation at $t^{\ell}$ is transformed into the map 
$U_0^{W(s^{\ell})}\longrightarrow \C_{s^{\ell}}$ obtained by evaluation at 
$s^{\ell}$. This proves the claim.
\end{proof}
\subsection{}
We now prove a crucial result which estimates the dimension of $D_t$. It will 
follow from the proof of Theorem \ref{mainqgthm} that this lower bound is precisely
the dimension of $D_t$.
\begin{lem}
For $t\in T$ the algebra $D_{t}$ has dimension at least $[W(t^{\ell}) : W(t)]$.
\end{lem}
\begin{proof}
For ease of notation let $R=\ULG$. Let $Q(R)$ be the quotient field of $R$ and 
suppose $G$ is a finite group acting faithfully on $R$ by automorphisms. Note 
that $Q(R^G) = Q(R)^G$. Indeed $Q(R^G)\subseteq Q(R)^G$. For the opposite 
inclusion observe that since $R$ is a finite module over $R^G$ we have that 
$Q(R) = R\otimes_{R^G} Q(R^G)$. Thus if $\alpha = xy^{-1}\in Q(R)^G$ we may 
assume without loss of generality that $y\in R^G$. It follows that $x\in R^G$ as 
required. Since $G$ acts faithfully on $R$ we have that
\begin{equation}
\label{eq:galois}
[Q(R):Q(R^G)]=[Q(R):Q(R)^G]= |G|,
\end{equation}
by \cite[Theorem 58]{rotman}.

Now consider the cases where $G=W(t)$ and $G=W(t^{\ell})$. Here we have
\[
[Q(R):Q(R^{W(t^{\ell})})] = [Q(R):Q(R^{W(t)})][Q(R^{W(t)}):Q(R^{W(t^{\ell})})].
\]
It follows from (\ref{eq:galois}) that
\[
[Q(R^{W(t)}):Q(R^{W(t^{\ell})})] = [W(t^{\ell}):W(t)].
\]
So generically, for $\mathfrak{m}$ a maximal ideal of $R^{W(t^{\ell})}$ the 
algebra $R^{W({t})}\otimes_{R^{W(t^{\ell})}} \mathbb{C}_{\mathfrak{m}}$ has dimension 
$[W(t^{\ell}):W(t)]$. Moreover, since the dimension function is upper 
semicontinuous and $R$ is an integral domain, in the non-generic case the 
dimension of $R^{W({t})}\otimes_{R^{W(t^{\ell})}} \mathbb{C}_{\mathfrak{m}}$ is at least   $[W(t^{\ell}):W(t)]$ as required.
\end{proof}

\subsection{}
\label{mainqgthm}
Let $\chi \in (U^-  \times U^+)\rtimes T$, so that $\mathfrak{m}_{\chi} \in \text{Maxspec}(Z_0)$. Define $U_{\chi} = \UE/\mathfrak{m}_{\chi}\UE$ and $Z_{\chi} = Z/\mathfrak{m}_{\chi}Z$. Much more easily than in the modular setting, we prove:
\begin{lem}
$Z_{\chi}$ is a subalgebra of $Z(U_{\chi})$.
\end{lem}
\begin{proof}
As $Z$-modules, $Z$ is a direct summand of $\UE$ since $Z$ is the image of the reduced trace map and we are in characteristic zero. Thus $\mathfrak{m}_{\chi}\UE \cap Z = \mathfrak{m}_{\chi}Z$, as required.
\end{proof}
\indent Continue with an element ${\chi}$ of $(U^-\times U^+) \rtimes T $. By Theorem \ref{notqg}.2, we can assume, without loss 
of generality, that $\chi \in B^-$. Then we have a unique decomposition $\chi = 
\chi_s \chi_u$ where $\chi_s \in T$ and $\chi_u \in U^-$. Consider $\{ \lambda \in T : \lambda^{\ell}\in W\chi_s^2 \}$. This subset of $T$ has a natural $W$-action and we will let $R_{\chi}$ denote a set of 
representatives for the $W$-orbits.

\begin{thm}
Keep the above notation, so in particular $\chi = \chi_s\chi_u \in B^-$. Then $Z_{\chi}$ is isomorphic to $\oplus_{\lambda \in R_{\chi}}D_{\lambda}$.
\end{thm}
\begin{proof}
By definition we have $Z_{\chi} = Z(\UE) \otimes_{Z_0} \mathbb{C}_{\chi}$. By Theorem 
\ref{notqg}.3
\begin{equation}
\label{zedchi}
 Z_{\chi}  \cong Z_1\otimes_{Z_1\cap Z_0} \mathbb{C}_{\chi},
\end{equation} 
where $Z_1\cap Z_0 \longrightarrow \mathbb{C}_{\chi}$ is obtained by restriction of 
$\chi$. Since $Z_1 \subset \UO + \UE\mathfrak{n}^+$ and since $\chi \in B^-$ it follows 
that $Z_1\cap Z_0 \longrightarrow \mathbb{C}_{\chi}$ depends only $\mathbb{C}_{\chi_s}$ and that 
this map can be evaluated by considering only the $U_0$ component of elements of $Z_1\cap 
Z_0$ written in triangular form. Since $\ULG$ is fixed pointwise by the map $\gamma$ which takes 
$K_{\lambda}$ to $\ep^{-(\rho,\lambda)}K_{\lambda}$, the last part of Theorem \ref{notqg} implies that the Harish-Chandra map 
$\psi$, of (\ref{eq:HCmap}) yields an isomorphism of $Z_1\otimes_{Z_1\cap Z_0} \mathbb{C}_{\chi}$ with $\UGW\otimes_{\ULGW} \mathbb{C}_{\chi_s}$ where $\ULGW \longrightarrow 
\mathbb{C}_{\chi_s}$ is restriction.

Consider the commutative diagram
\[
\begin{CD}
T@>(-)^{2} >> T @>a >> T/W \\
@VV (-)^{\ell} V @VbV(-)^{\ell}V @VcV(-)^{\ell}V \\
T@>(-)^{2} >> T @>d>> T/W
\end{CD}
\]
induced from the inclusions
\[
\begin{CD}
\UO @<<< \UOG @<<<  \UGW \\
@AAA @AAA @AAA \\
\UL @<<< \ULG @<<< \ULGW . \\
\end{CD}
\]
Then, by the above, $Z_{\chi} \cong c^*(d(\chi_s^2))$. The closed points are 
given by
\[
c^{-1}d(\chi_s^2) = ab^{-1}d^{-1}d(\chi_s^2) = \{ t\in T : t^{\ell} \in 
W\chi_s^2 \} / W,
\]
so $R_{\chi}$ parametrises the components of $c^*(d(\chi_s^2))$.

Let $t\in R_{\chi}$ and let $V_t$ be the component of $c^*(d(\chi_s^2))$ 
containing the orbit associated to $t$. Now consider the diagram 
\[
\xymatrix{
T \ar[d]_b \ar[r]^(.3){e} & T/W(t) 
\ar[d]_s \ar[rr]^f &&  T/W \ar[d]_c 
\\
T \ar[r] & T/W(t) 
\ar[r]^m & T/W(t^{\ell}) \ar[r]^o & T/W .
}
\]
Since $W(t)$ is the stabiliser of $t$, the map $f$ is a covering near $e(t)$ so 
the component of $(cf)^*(d(\chi_s^2))$ containing $e(t)$ is isomorphic to $V_t$.

Consider next $(cf)^*(d(\chi_s^2)) = (oms)^*(d(\chi_s^2))$. Since $W(t^{\ell})$ 
is the stabiliser of $mse(t)$ the map $o$ is a covering near $mse(t)$ so the 
component of $o^*(d(\chi_s)^2)$ containing $mse(t)$ is precisely the (reduced) 
point $mse(t)$.

Now $m^*o^*(d(\chi_s)^2)$ has a unique closed point since $W(t)$ stabilises 
$se(t)$. The neighbourhood of this point is given by the kernel of the map
\[
U_{0,\ell}^{W(t)\ltimes \Gamma} \longrightarrow U_{0,\ell}^{W(t)\ltimes \Gamma} 
\otimes_{U_{0,\ell}^{W(t^{\ell})\ltimes \Gamma}}\mathbb{C}_{mse(t)}.
\]
The tensor product on the right is isomorphic to $D_{t}$.

Now $V_t$ is a neighbourhood of $e(t)$ and we have shown that $D_t$ is a 
subalgebra of $\mathcal{O}(V_t)$. Therefore $\oplus_{t\in R_{\chi}}D_t$ is a 
subalgebra of $\oplus_{t\in R_{\chi}}\mathcal{O}(V_t) \cong Z_{\chi}$. Noting 
that the number of elements in the $W$-orbit associated to $t$ such that 
$t^{\ell} = \chi_s^2$ is precisely $[W(t^{\ell}):W(t)]$ we see from Lemma 
\ref{size} that
\[
\dim (\oplus_{t\in R_{\chi}}D_t) \geq \sum_{t\in R_{\chi}} [W(t^{\ell}):W(t)] = 
|\{ t\in T: t^{\ell} \in \chi_s^2 \}| = \ell ^r.
\]
Since $\dim (Z_\chi) = \ell^r$ by (\ref{zedchi}) and Theorem \ref{notqg}.5, the theorem is proved.
\end{proof}
\subsection{Regular characters}
\label{regqg}
Using Theorem \ref{mainqgthm} we can completely describe the representation 
theory of the algebras $U_{\chi}$ for $\chi$ a regular character. By definition 
the regular characters $\chi$ correspond to the points of $\text{Maxspec}(Z_0)$ 
lying in symplectic leaves of maximal dimension. It is shown in \cite[Theorem 
24.1]{decpro3} that the irreducible modules of $U_{\chi}$ for $\chi$ regular all 
have dimension equal to $\ell ^N$, the PI degree of $\UE$. Therefore we can 
apply Corollary \ref{regalgs} to deduce the first part of the following corollary. Part 3    is implicit in \cite{decpro3}.
\begin{cor}
1. Let $\chi$ be a regular character of $\mathrm{Maxspec}(Z_0)$. Then there is an algebra 
isomorphism
\[
U_{\chi} \cong \mathrm{Mat}_{\ell^N}(Z_{\chi}).
\]
2. If $\chi \in U^-$ is regular there is an isomorphism
\[
U_{\chi} \cong \bigoplus_{t} \mathrm{Mat}_{\ell^N}\left( U_0^{W(t)}\otimes 
_{U_0^{W}} \C_1\right) ,
\]
where $t$ runs over a set of orbit representatives for the action of $W$ on the 
set $\{ t\in T: t^{\ell} = 1\}$.\\
3. If $\chi\in T$ is regular then $U_{\chi}$ is semisimple, namely 
\[U_{\chi} \cong \oplus^{\ell ^r} 
\mathrm{Mat}_{\ell^N}( \C ).
\]  
\end{cor}
\begin{proof}
3. Suppose that $\chi\in T$  is regular. Then by definition (see Theorem \ref{notqg}) $\pi(\chi) = \chi^2$ is a regular semisimple 
element of $G^0\subset G$ and so $W(\chi^2) = 1$. Therefore, for $t\in T$ such 
that $t^{\ell} = \chi^2$, we have $D_t = \C$, and no two such $t$ can be $W$-conjugate, so there are $\ell^r$ summands in the decomposition of $Z_{\chi}$ given by Theorem \ref{mainqgthm}.  
\end{proof}
\begin{rem}
The above corollary also shows that the algebras $D_t$ are
Frobenius. To see this one need only note that given $t\in T$ there
exists $u\in U^-$ such that $g=tu$ is the Jordan decomposition of a
regular element of $G$. Now follow the proof of Theorem \ref{mirrum}.2.
\end{rem}
\subsection{}
It's instructive to do a simple calculation for the case of 
$U_{\ep}(\mathfrak{sl}_2)$ and $\chi_s = 1$. Then the elements $t\in T$ such 
that $t^{\ell} = 1$ are parametrised by the integers from $0$ to $\ell -1$. 
Specifically $i$ corresponds to $t_i =\bigl( \begin{smallmatrix} \ep^i & 0 \\ 0& 
\ep^{-i} \end{smallmatrix} \bigr)$ for $0\leq i \leq \ell -1$. Then we have 
$W(t_i^{\ell}) = W$ whilst
\begin{equation*}
W(t_i) = 
\begin{cases}
W \quad & \text{if }i= 0, \\
1 \quad & \text{if }i\neq 0.
\end{cases}
\end{equation*}
Observe that $U_0^W = \C[ K_{\varpi} + K_{\varpi}^{-1}]$ and that the kernel of 
$U_0^W \longrightarrow \C_1$ is the ideal generated by $K_{\varpi} + 
K_{\varpi}^{-1} - 2$. We deduce that
\[
D_t \cong 
\begin{cases}
\C \quad & \text{if } i=0,\\
\frac{\C[K_{\varpi}^{\pm 1}]}{<(K_{\varpi}-1)^2>} \quad & \text{if }i\neq 0.
\end{cases}
\]
\indent In particular if $\chi$ is regular unipotent we see from Corollary (\ref{regalgs}) that 
\[
U_{\chi}(\mathfrak{sl}_2)\quad \cong \quad \mathrm{Mat}_{\ell}(\C) \oplus 
\bigoplus^{\frac{\ell-1}{2}}\mathrm{Mat}_{\ell}\left( \frac{\C[X]}{<X^2>}\right)	.
\]
\subsection{Blocks and baby Verma modules}
\label{quanbab}
We
shall briefly discuss baby Verma modules for $\UE$. Using Theorem \ref{notqg}
we can concentrate on characters $\chi = \chi_u\chi_s\in B^-$ without
loss of generality. Let $t\in T$ be an $\ell^{\text{th}}$ root of
$\chi_s$. The baby Verma module $V_{\chi}(t)$ is defined to be the
induced module $U_{\chi}\otimes_{U_{\chi}^{\geq 0}} \mathbb{C}_t$
where $\mathbb{C}_t$ is the one-dimensional $U_{\chi}^{\geq 0}$-module
with $E_{\alpha}.1 = 0$ and $K_{\lambda}.1 = \lambda (t)$. We will
write $v_t$ for the element $1\otimes 1\in V_{\chi}(t)$.

Recall from (\ref{notqg}) the map $\gamma :U_0 \longrightarrow U_0$ which sends
$K_{\lambda}$ to $\ep^{-(\rho,\lambda)}K_{\lambda}$ and let
$\gamma^*:T\longrightarrow T$ be the corresponding morphism. Recall also the Harish-Chandra map $\psi = \gamma\circ \pi_{U_{0}}$. For ease
of notation for $w\in W$ and $t\in T$ we'll write
\[
w{\scriptscriptstyle \bullet} t = \gamma^* \circ w \circ (\gamma^*)^{-1}(t).
\]
Since
$Z_1\subseteq U_0 + \sum_{\alpha\in \Delta}\UE E_{\alpha}$ we see that
$z\in Z_1$ annihilates $V_{\chi}(t)$ if and only if $\pi_{U_0}(z)(t)=
0$ if and only if $\psi(z)((\gamma^*)^{-1}(t))=0$. In particular
$\text{Ann}_{Z_1}(V_{\chi}(t))$ is a maximal ideal of $Z_1$ and, since
$\psi(Z_1)\subseteq U_0^W$, we deduce
\begin{equation}
\label{eq:quantlinkage}
\textit{if $V_{\chi}(t)$ and
$V_{\chi}(u)$ belong to the same block then $t=w{\scriptscriptstyle \bullet} u$ for some
$w\in W$.}
\end{equation}
Moreover, every
maximal ideal of $Z$ lying over $\mathfrak{m}_{\chi}$ occurs as the
annihilator of a baby Verma module, since (as in the modular setting)
we can find a baby Verma module mapping onto any irreducible
$U_{\chi}$-module $V$, by viewing $V$ as a $U_{\chi}^{\geq 0}$-module,
so finding a highest weight vector in $V$, and taking the baby Verma
module which the latter generates. In the light of
(\ref{eq:quantlinkage}) and noting that irreducibles in the same block have the same central annihilator we can now record the analogue for $\UE$ of Theorem \ref{gblocks}. Recall the definition of $R_{\chi}$ from (\ref{mainqgthm}).
\begin{thm}
Let $\chi = \chi_s \chi_u \in B^-$. Then the blocks of $U_{\chi}$ are
in natural bijection with $R_{\chi}$.
\end{thm}
We take this  opportunity to point
out that the claim made in \cite[(20.3)]{decpro3} that baby Verma
modules are indecomposable is false, as can be seen by considering the
following example, analogous with the modular situation
\cite[(6.9)]{jan6}. Let $\g = \mathfrak{sl}_3$ and assume that $\ell$ is prime to $2$ and
$3$. Consider the unipotent element of $SL(3)$
\[
u= \left( 
\matrix
1 & 0 & 0 \\
0 & 1 & 0 \\
1 & 0 & 1 
\endmatrix
\right).
\]
It is easy to check that this is a subregular element - that is,
$\dim (SL(3).u) = 4$. Any character, $\chi \in \mathrm{Maxspec}(Z_{0})$, lying over $u$ must
satisfy the equalities
\[
\chi (F_{\alpha_1}^{\ell}) = 0 = \chi (F_{\alpha_2}^{\ell}), \quad
\chi(K_{\lambda}^{2\ell}) = 1, \quad \chi (E_{\alpha_1}) =\chi
(E_{\alpha_2}) = \chi(E_{\alpha_1+\alpha_2})=0.
\] 
Let $n$ be the inverse of $3$ in $\mathbb{Z}/\ell \mathbb{Z}$ and let $t \in \mathrm{Maxspec}( U_0)$ be such that $t(K_{\varpi_1}) =
\ep^{n-2}$ and $t(K_{\varpi_2}) =
\ep^{2n-2}$. Then, following \cite[(1.3.2)]{deckac5} there
are maps
\[
\theta_1:V_{\chi}(t_1) \longrightarrow V_{\chi}(t),
\]
and 
\[
\theta_2:V_{\chi}(t_2) \longrightarrow V_{\chi}(t)
\]
given by
\[
v_{t_1} \longmapsto
F_{\alpha_1}^{\ell - 1}v_{t},
\]
and
\[
v_{t_2} \longmapsto
F_{\alpha_2}v_{t},
\]
where $t_1(K_{\varpi_1}) = \ep^{n-3}, t_1(K_{\varpi_2})=\ep^{2n-2}$
and $t_2(K_{\varpi_1}) = \ep^{n-2}, t_2(K_{\varpi_2})= 1$. It is clear
that the set  
\begin{equation}
\label{eq:furry}
\{  F_{\alpha_1}^{i}F_{\alpha_1+\alpha_2}^{j}F_{\alpha_2}^kv_{t}:
0\leq i,j\leq \ell -1 , 1 \leq k\leq \ell-1 \}
\end{equation}
is a basis of $\im (\theta_2)$, and similarly for
$\im (\theta_1)$. Therefore $\dim \im
(\theta_1) = \ell^2$ whilst $\dim \im (\theta_2) =
\ell^2 (\ell - 1)$. By \cite{cant} all simple $U_{\chi}$-modules have dimension divisible by
$\ell^2$, so $\im (\theta_1)$ is simple. Thus $\im (\theta_1)$ and $\im (\theta_2)$ intersect non-trivially if and only
if $F_{\alpha_1}^{\ell - 1}v_{t} \in \im (\theta_2)$. Therefore we
deduce from
(\ref{eq:furry}) that we have a decomposition 
\[
V_{\chi}(t) = \im (\theta_1) \oplus \im (\theta_2) .
\]
\subsection{}
\label{projqg}
We consider now the quantum analogue of Corollary \ref{babyproj}.
\begin{thm}
Let $\chi \in B^-$ and continue with the notation introduced in (\ref{quanbab}). Then the following statements are equivalent:\\
(i) $V_{\chi}(t)$ is  a projective $U_{\chi}$-module.\\
(ii) $V_{\chi}(t)$ is a projective and irreducible $U_{\chi}$-module.\\
(iii) $\mathrm{Ann}_Z(V_{\chi}(t))$ is an unramified point of $\mathrm{Maxspec}(Z)$.\\
(iv) for every root $\alpha\in\tilde{\Delta} = \Delta\cup\{- \alpha_0\}$, 
$\alpha(t)^{2\ell} = 1$ implies $\alpha (t)^2 = \ep^{-(2\rho,\alpha)}$.
\end{thm}
\begin{proof}
(ii)$\Longrightarrow$(i): Trivial.\\
(i)$\Longrightarrow$(ii):
Suppose that $V_{\chi}(t)$ is projective. Then its restriction to $U_{\chi}^{+}$ is also 
projective. Since $\chi\in B^-$, $U_{\chi}^{+}$ is a local ring 
of dimension $\ell ^N = \dim (V_{\chi}(t))$. Thus, as $U_{\chi}^{+}$-module $V_{\chi}(t)$ is free of rank one. Moreover, being scalar local and self-injective, $U_{\chi}^{+}$ has a 
simple socle which corresponds to $1\otimes \C_t$. Therefore any non-zero 
$U_{\chi}$-submodule of $V_{\chi}(t)$ must contain $1\otimes \C_t$, proving that 
$V_{\chi}(t)$ is irreducible.\\
(ii)$\Longleftrightarrow$(iii): Theorem (\ref{azunram}) applies in view of Theorem (\ref{az=sing})2.\\
(iii)$\Longleftrightarrow$(iv) In terms of 
$\UGW$ the central character of $V_{\chi}(t)$ is described by the Harish-Chandra 
map, $\psi$, that is $\UGW$ acts on $V_{\chi}(t)$ by restriction from the action of $U_{0}$ given by  
$K_{\mu}\cdot 1 = \ep^{-(\rho, \mu)}\mu(t)$. Call the torus element 
corresponding to this character $u$. Considering this as an element of $T\cong 
\text{Spec}(U_0^{\Gamma})$ we consider $u^2$. By Theorem \ref{mainqgthm} the 
central character of $V_{\chi}(t)$ is unramified if and only if $W(u^{2\ell}) = 
W(u^2)$. By (\ref{hearmysong}) we can assume without loss of generality that $W(u^{2\ell})$ is 
generated by reflections coming from a proper subset of $\tilde{\Delta}$. Now 
$s_{\alpha}(u^{2\ell}) = u^{2\ell}$ if and only if $K_{s_{\alpha}\mu}(u^{2\ell}) 
= K_{\mu}(u^{2\ell})$ for all $\mu \in X$. This occurs if and only if 
$K_{\alpha}^{< \mu, \check{\alpha} >}(u^{2\ell}) = 1$ for all $\mu\in X$. 
Similarly for $W(u^2)$. Therefore we deduce that $W(u^{2\ell}) = W(u^2)$ if and 
only if $K_{\alpha}(u)^{2\ell} = 1$ implies $K_{\alpha}(u)^2 = 1$ for all roots 
$\alpha\in \tilde{\Delta}$. By the definition of $u$ this last condition is 
equivalent to $\alpha (t)^{2\ell}=1$ implies $\alpha (t)^2 = \ep^{-(2\rho , 
\alpha)}$ for all $\alpha\in \tilde{\Delta}$.
\end{proof}
\subsection{}
\label{exceptional}
In the next few sections we will give a necessary condition for the
simplicity of a baby Verma module. To begin with we must recall the notion of exceptional
elements of $G$. A semisimple element $g\in G$ is called
\emph{exceptional} if the centraliser $C_G(g)$ of $g$ in $G$ has a
finite centre. An arbitrary element $g\in G$ is called
\emph{exceptional} if its semisimple part, with respect to the Jordan
decomposition, is exceptional. Finally, a $U_{\chi}$-module is called
\emph{exceptional} if $\chi \in \mathrm{Maxspec}(
Z_0)$ has $\pi(\chi)$ exceptional (where $\pi$ is the covering map of Theorem \ref{notqg}).
\begin{lem}
Let $\alpha_0 = \sum_{i=1}^r a_i\alpha_i$ be the highest root in
$\Phi^+$.\\
1.  For $m = 1, \ldots , r$, define elements $h_m \in \h$ by $\alpha_j(h_m) =
\delta_{jm}$. Then the elements
\begin{equation}
\label{eq:exp}
s_m = \mbox{exp}(2\pi i h_m/a_m) \in T\subset G
\end{equation}
and $s_0 = 1$ are exceptional semisimple elements and any exceptional
semisimple element is conjugate to one of the $s_m$.
\newline
\noindent
2.  Fix $1\leq m \leq r$ and let $\beta_m$  be the minimal root
in $\Phi^+$ whose $\alpha_m$-coefficient is $a_m$. Then the
centraliser of $s_m$ is the reductive group generated by $T$ and the
root subgroups $U_{\alpha_i}$, $i\neq m$, and $U_{\beta_m}$.
\end{lem}
\begin{proof}
Part 1 can be deduced from \cite[Chapter 8]{kacILA}. Since $G$ is simply-connected it follows from
\cite[Theorem 2.2 and Theorem 2.11]{humCC} that 
$C_G(s_m)$ is the reductive group generated by $T$ and the root
subgroups $U_{\beta}$ where $\beta(s_m)=1$. If
$\beta = \sum_{i=1}^r b_i\alpha_i$ it
follows from (\ref{eq:exp}) that $\beta(s_m) = 1$ if and only if $a_m$ divides $b_m$.
As $a_m$ is a coefficient of the highest root of $\Phi$ 
this occurs if and only if $b_m = a_m$ or $b_m=0$.
\end{proof}
\indent Part 1 of the lemma shows that a complete set of
representatives of the conjugacy classes of exceptional elements is
given by $\{ s_m u \}$ with $s_m$ as above and $u\in C_G(s_m)$ a
representative of a conjugacy class of unipotent elements in the
centraliser of $s_m$. 
\subsection{}
\label{quanred}
Fix $\chi=\chi_u\chi_s\in \mathrm{Maxspec} (Z_0)$ lying over $g=g_ug_s\in G$. Let
\[
\h' = \text{Lie}(Z(C_G(g_s))) \subseteq \h,
\]
and let
\begin{equation}
\label{eq:barbaracartland}
\Phi' = \{ \beta \in \Phi : \beta \text{ vanishes on } \h '\} =
Q ' \cap \Phi,
\end{equation}
where $\Delta' = \Phi' \cap \Delta$ and $Q'= \mathbb{Z}\Delta'$.   

Let $U_{\ep}'(\g)$ be the subalgebra of $\UE$ generated by
$U_0$ and the elements
$E_i$ and $F_i$ such that $\alpha_i\in \Delta'$ and let $U_{\chi}'$ be
the corresponding subalgebra of $U_{\chi}$. Similarly let
$\tilde{U}_{\ep}'(\g)$ be the subalgebra generated by $U^+ \otimes U_0  $
and the elements $F_i$ such that $\alpha_i\in \Delta'$ and let
$\tilde{U}_{\chi}'$ be the corresponding subalgebra of $U_{\chi}$. Finally let
$U_{\ep}''(\g)$ be the subalgebra of $\UE$ generated by the elements
$1$ and $K_{\alpha_i}$, $E_i$ and $F_i$ such that $\alpha_i\in
\Delta'$ and let 
$U_{\chi}''$ the corresponding subalgebra of $U_{\chi}$.
\begin{thm}\cite[Theorem 8]{deckac100}
Let $V$ be a simple $U_{\chi}$-module.
\newline
\noindent
1.  The $U_{\chi}$-module $V$ contains a unique simple
$\tilde{U}_{\chi}'$-submodule $V'$. On restriction, $V'$ is a simple
$U_{\chi}'$-module.
\newline
\noindent
2.  There is an isomorphism of $U_{\chi}$-modules
\[
V \cong U_{\chi}\otimes_{\tilde{U}_{\chi}'} V'.
\]
\noindent
3.  The map $V\longmapsto V'$ is a bijection between the simple
$U_{\chi}$-modules and the simple $U_{\chi}'$-modules.
\newline
\noindent
4.  As a $U_{\chi}''$-module $V'$ is an exceptional simple module.
\end{thm} 
\subsection{}
\label{simpqg}
Suppose $\chi = \chi_u\chi_s \in \mathrm{Maxspec} (Z_0)$ lies over $g=g_ug_s$. After conjugating we can and shall assume that $g_s \in T$. As
$C_G(g_s)$ is reductive \cite[Theorem 2.2]{humCC}, we can write it as a product
\[
C_G(g_s) = Z \cdot G_1\cdots G_n,
\]
where $Z$ is the centre of $C_G(g_s)$ and each $G_i$ is a simple algebraic group. Let $\Phi_i$ be the root
system of $G_i$ and let $\Delta_i$ be a set of simple roots of
$\Phi_i$ lying in $\Phi^+$. Write $g_u\in C_G(g_s)$ as a product of
unipotent elements $g_i$ in $G_i$,
\[
g_u = g_1\ldots g_n.
\]
\begin{thm}
Keep the above notation. If $V_{\chi}(t)$ is a simple
$U_{\chi}$-module then for each $i$, $1\leq i\leq n$, either $g_i$ is
regular in $G_i$ or $\alpha(t)^2 = \ep^{-2(\rho , \alpha)}$ for all
$\alpha \in \Delta_i$.
\end{thm}
The proof of this theorem is quite technical and will occupy the next
four sections. As we will see this theorem is enough for us to
describe the fully Azumaya points of $\mathrm{Maxspec} Z_0$.

\begin{rem}
Using Theorem \ref{az=sing}(2) this result can be considered as a
step towards describing the smooth locus of
$\text{Maxspec}(Z)$. In order to complete the description we would
have decide whether the above condition implies that $V_{\chi}(t)$
is simple. It may be possible to a certain extent to imitate the
analysis of the
Lie algebra case in \cite{kry}. This approach shows
that an element $(g_ug_s,tW)\in G\times_{T/W}T/W$ is smooth if and only
if it is smooth in $C_G(g_s)\times_{T/W}T/W$. The problem now arises,
however, that $C_G(g_s)$ is not in general simply-connected and so the
analysis of the map $C_G(g_s)\longrightarrow T/W$ in \cite[Section
8]{ste} is unavailable. Moreover the centre of an adjoint quantised
enveloping algebra is unknown in general, \cite[Remark 6.4]{deckacpro4}. One might hope, however, that if $\ell$ is
prime to the index of connection of $\g$ the necessary condition given in Theorem \ref{simpqg} is also sufficient.
\end{rem}
\subsection{}
\label{reduce}
We begin by reducing to the exceptional case. We will use the notation
of Section \ref{quanred}. In addition, let $\g'$ be the
(semisimple) Lie algebra associated to the root lattice $\Phi'$. Let
$X'$ denote the weight lattice of $\g'$ and recall that $Q'$ is the
root lattice. 

Let $H$ be the reductive algebraic subgroup of $G$ generated by $T$
and the root subgroups $U_{\alpha}$ for $\alpha\in \Phi'$. Let $G'$
be the simply-connected algebraic group associated with $\g'$, and let
$T'$ be a maximal torus of $G'$. There is a projection map $p$,
\[
\xymatrix{
G \supseteq 
H \ar@{->>}[r]^(0.6){p} & G'
},
\]
which sends $T$ to $T'$ and whose kernel consists of central elements of $H$. In
particular $p$ is bijective on $U^-\cap H$. Let
${U'}^-=p(U^-\cap H)$. Since $C_G(g_s)\subseteq
H$, $p$ restricts to provide a surjective morphism
\begin{equation}
\label{eq:indianajones}
p: C_G(g_s) \twoheadrightarrow C_{G'}(p(g_s)).
\end{equation}
Since $p$ is injective on unipotent parts it follows that
$C_{G'}(g_s) = Z'\cdot p(G_1)\cdots p(G_n)$ and that $g_i$ is
regular in $G_i$ if and only if $p(g_i)$ is regular in $p(G_i)$. 

Let $\iota :U_{\ep}^0(\g') \cong \C [X'] \subseteq \C[X] \cong U_0$ be
the natural inclusion and let $\tilde{p}:T\longrightarrow T'$ be the
corresponding morphism. By construction this morphism is the
restriction of $p$ to $T$. It follows from (\ref{eq:indianajones})
and the fact that the elements in the kernel of $\tilde{p}$ are
central in $C_G(g_s)$ that $\alpha(t)=\alpha(\tilde{p}(t))$ for all
$t\in T$ and $\alpha$ belonging to the roots of $C_G(g_s)$. 
  
In summary we have shown for all $i$, $1\leq i\leq n$,
\begin{equation}
\label{eq:anniehall}
\textit{$g_i$ is regular in $G_i$ if and only if $p(g_i)$ is regular in
$p(G_i)$}
\end{equation}
and
\begin{equation}
\label{eq:caddyshack}
\textit{$\alpha(t)^2=\ep^{-2(\rho,\alpha)}$ for all $\alpha\in \Delta_i$
if and only if $\alpha(\tilde{p}(t))^2=\ep^{-2(\rho,\alpha)}$ for all $\alpha \in
\Delta_i$.}
\end{equation}

There are inclusions of algebras $\UE \supseteq U_{\ep}'(\g) \supseteq
U_{\ep}''(\g)$. Moreover there is an isomorphism $U_{\ep}''(\g) \cong
U_{\ep , Q'}(\g')$. (To be completely precise it may be necessary to
replace $\ep$ by $\ep^{d_i}$ for some $i$, $1\leq i\leq r$. This, however, will
have no bearing on the reduction since $d_i$ is prime to $\ell$.)
Thus there is an injection 
$U_{\ep}''(\g)\longrightarrow U_{\ep}(\g')$ into the simply-connected form of the quantised enveloping algebra which is bijective on the
positive and negative graded pieces of $U_{\ep}''(\g)$. Let
$Z_0^-(\g') \subseteq Z_0^-$ and $Z_0^0(\g')\subseteq Z_0^0$ be the
natural inclusions. These yield a map
\[
\begin{CD}
\mathrm{Maxspec} Z_0^-\times \mathrm{Maxspec} Z_0^0 @>\theta >> \mathrm{Maxspec} Z_0^-(\g ')\times \mathrm{Maxspec}
Z_0^0(\g').
\end{CD}
\]
Let $\pi': \text{Maxspec}Z_0^-(\g ')\times \mathrm{Maxspec}Z_0^0(\g')
\longrightarrow {U'}^-\times T'$ be the analogue of the map in Theorem
\ref{notqg} for $U_{\ep}(\g')$. Then by construction we have $p(\pi
(\chi_s)) = \pi ' (\theta(\chi_s))$. On the other hand if we choose a
reduced expression of the longest word in $W$ which begins with a
reduced expression for the longest word in the Weyl group of $\g'$
then we see that
\begin{eqnarray*}
\pi'(\theta(\chi_u)) &=& \text{exp}(\theta (\chi_u)(y_{\beta_{N'}}f_{\beta_{N'}}))\ldots
\text{exp}(\theta (\chi_u)(y_{\beta_1}f_{\beta_1}))\\
& =&
\text{exp}(\chi_u(y_{\beta_{N'}})f_{\beta_{N'}})\ldots
\text{exp}(\chi_u(y_{\beta_1})f_{\beta_1}),
\end{eqnarray*}
whilst
\begin{equation}
\label{eq:piggy}
p(\pi(\chi_u)) = p(\text{exp}(\chi_u(y_{\beta_N})f_{\beta_N})\ldots
\text{exp}(\chi_u(y_{\beta_1})f_{\beta_1})).
\end{equation}
Since $\chi_s(K_{\beta}^{2\ell})\neq 1$ for all $\beta\notin \Phi'$,
Lemma \ref{hearmysong} together with the definition of $y_{\beta}$ shows that $\chi_u(F_{\beta}^{\ell}) = 0$ for
all such $\beta$. Therefore (\ref{eq:piggy}) becomes
\begin{eqnarray*}
p(\pi(\chi_u))& =& p(\text{exp}(\chi_u(y_{\beta_{N'}})f_{\beta_{N'}})\ldots
\text{exp}(\chi_u(y_{\beta_1})f_{\beta_1})) \\
&= &\text{exp}(\chi_u(y_{\beta_{N'}})f_{\beta_{N'}})\ldots
\text{exp}(\chi_u(y_{\beta_1})f_{\beta_1}),
\end{eqnarray*}
showing that $p(\pi(\chi_u)) = \pi'(\theta(\chi_u))$.

Now suppose $V=V_{\chi}(t)$ is
a simple $U_{\chi}$-module. By Theorem \ref{quanred}.2 there exists a
simple $\tilde{U}_{\chi}'$-module, written $V'$, such that there is an
isomorphism
\begin{equation}
\label{eq:fishliveinconcert}
V\cong
U_{\chi}\otimes_{\tilde{U}_{\chi}'} V',
\end{equation}
showing in particular that $V'$ has dimension $\ell^m$ where $2m=|\Phi'|$. 
Moreover, by Theorem
\ref{quanred}.4, $V'$ restricts to a simple $U_{\chi}''$-module. It follows
from the discussion above that $U_{\chi}''$ is isomorphic to the adjoint subalgebra
of $U_{\theta(\chi)}(\g')$, generated by $K_{\alpha}^{\pm}$,
$E_{\alpha}$ and $F_{\alpha}$ for $\alpha\in Q'$. Consider the baby
Verma $U_{\ep}(\g')$-module $V_{\theta(\chi)}(\tilde{p}(t))$. By
weight considerations its
restriction to $U_{\chi}''$ maps to $V'$, so,
since $V'$ is a simple $U_{\chi}''$-module of dimension $\ell^m$, this
map must be an isomorphism.  Therefore
$V_{\theta(\chi)}(\tilde{p}(t))$ must also be simple as a $U_{\theta(\chi)}(\g')$-module.
Since $\theta(\chi)$ lies over
$p(g_u)p(g_s)$ the reduction now follows from (\ref{eq:anniehall}) and (\ref{eq:caddyshack}).
\subsection{}
\label{calculation}
We begin with a combinatorial lemma. Recall that, for $1\leq m\leq r$, $\beta_m$ is the
minimal root in $\Phi^+$ whose $\alpha_m$-coefficient equals that of
the $\alpha_m$ coefficient in $\alpha_0$.
\begin{lem}
For $m$, $1\leq m\leq r$, there exists $w_m\in W$ such that
\newline
(i) $\beta_m= w_m (\alpha)$ for some $\alpha \in \Delta$,
\newline
(ii) each element $\gamma \in \Phi^+\cap w_m(\Phi^-)$ has 
the form $\gamma = \sum_{i=1}^r c_i\alpha_i$ where $c_m \gneq 0$,
\newline
(iii) each element $\gamma \in \Phi^+\cap w_m(\Phi^-)$ is strictly
less than $\beta_m$ in the usual partial ordering on $\Phi^+$.
\end{lem}
\begin{proof}
Let $w\in W$ have reduced expression $w= s_{i_1}\ldots s_{i_t}$ and
let $\gamma_1 = \alpha_{i_1}$ and $\gamma_k = s_{i_1}\ldots s_{i_{k-1}}(\alpha_{i_k})$ for $1 < k
\leq t$. Then, by \cite[Chapter 1, Proposition 3.6]{hilCG}, 
\[
\Phi^+ \cap w(\Phi^-) = \{ \gamma_k : 1\leq k \leq t \}.
\]
The result now follows from case-by-case analysis using the explicit
elements $w_m\in W$ given in the table in
the appendix.
\end{proof}
\subsection{}
\label{twist}
Let $w=s_{i_1}\ldots s_{i_t}\in W$ be a reduced expression and define $w_j = s_{i_1}\ldots
s_{i_j}$ for $1\leq j\leq t$. Let $\gamma_j =w_{j-1}(\alpha_{i_j})$
(where $\gamma_1 = \alpha_{i_1}$). Suppose  
that $\chi=\chi_u\chi_s \in \mathrm{Maxspec} (Z_0)$ is zero on elements of $Z_0^-$ of weight $-\ell\gamma_j$ for
$1\leq j\leq t$. 
Let $I$ be the ideal of $\UE$ defining $U_{\chi}$, that is the ideal
generated by the elements $E_{\beta}^{\ell},
F_{\beta}^{\ell} - \chi_u(F_{\beta}^{\ell})$ for $\beta \in \Phi^+$, and by
$K_{\lambda}^{\ell}-\chi_s(K_{\lambda}^{\ell})$ for $\lambda \in X$. Let $I^{\geq
0}\subseteq I$ be the ideal of $U^+ \otimes U^0$ generated by
$E_{\beta}^{\ell}$ and $K_{\lambda}^{\ell}-\chi_s(K_{\lambda}^{\ell})$. 
Note that
$T_{w_j}(E_{\beta}^{\ell}) \in I$ for all $\beta \in \Phi^+$ and $j=1, \ldots ,t$, where $T_{w_j}$ denotes the braid group automorphism of $\UE$ corresponding to $w_j$, \cite[Chapter 8]{janqg}. Indeed if  $w_j(\beta) \in \Phi^+$ this is
clear, whilst if $w_j(\beta)\in \Phi^-$ the element
$T_{w_j}(E_{\beta}^{\ell})$ is homogeneous of weight $-\ell\gamma_k$ for
some $1\leq k\leq j$, so, by hypothesis, lies in $I$ too. Thus there
is an inclusion 
\[
T_{w_j}(U_{\chi_{u}w_j^{-1}\chi_{s}}^{\geq 0}) \subseteq U_{\chi}(\g).
\]

For $w_j$ and $\chi$ as above and for $u$ an $\ell$th root of $w_j^{-1}\chi_{s}$, we can therefore define twisted baby Verma modules as
follows:
\[
V_{\chi}^{w_j}(u) = U_{\chi} \otimes_{T_{w_j}(U_{\chi_{u}w_j^{-1}\chi_{s}}^{\geq
0})} k_{u}.
\]
Observe that this module has the following universal property: if $V$ is a $U_{\chi}$-module and $0\neq v\in V$ is an element
such that $\C v \cong k_u$ as $T_{w_j}(U_{\chi_{u}w_j^{-1}\chi_{s}}^{\geq
0})$-modules then there is a $U_{\chi}$-homomorphism $V_{\chi}^{w_j}(u)
\longrightarrow V$ obtained by sending $1\otimes 1$ to $v$.
\begin{lem}
Assume $j\leq t-1$. Then there is a $U_{\chi}$-homomorphism 
\[
\theta_{j+1}:V_{\chi}^{w_j}(u) \longrightarrow
V_{\chi}^{w_js_{i_{j+1}}}(u'),
\]
where $u'(K_\lambda) =
\ep^{(\lambda,\gamma_{j+1})}u(K_{\lambda})$. The homomorphism is
induced by sending $1\otimes 1$ to
$T_{w_j}(E_{i_{j+1}}^{\ell-1})\otimes 1$.
\end{lem}
\begin{proof}
By the universal property it is enough to check that
$T_{w_j}(E_{i_{j+1}}^{\ell-1})\otimes 1\in V_{\chi}^{w_js_{i_{j+1}}}(u')$ is a vector of weight $u$ which is annihilated by $T_{w_j}(E_{\alpha_i})$ for all
$i$, $1\leq i\leq r$. The weight claim follows from the equality
\[
K_{\lambda}T_{w_j}(E_{i_{j+1}}^{\ell-1}) =
\ep^{(l-1)(\lambda,w_j\alpha_{i_{j+1}})}T_{w_j}(E_{i_{j+1}}^{\ell
-1})K_{\lambda}.
\]

It remains only to show that $T_{w_j}(E_k)$ annihilates
$T_{w_j}(E_{i_{j+1}}^{\ell-1})\otimes 1$ for $1\leq k\leq r$. If $k =
i_{j+1}$ this is clear so we can assume that $k\neq i_{j+1}$. We
pick a reduced expression for the longest word of $W$ starting with
the simple reflection $s_{i_{j+1}}$. Then, in $U_{\ep}^+$, there is a
commutation relation
\[
E_kE_{i_{j+1}}^{\ell -1} = \sum_{n=1}^{\ell -1} E_{i_{j+1}}^n u_n,
\]
where $u_n = T_{s_{i_{j+1}}}(u_n')$ for some $u_n'\in U^+$ of weight
$s_{\alpha_{i_{j+1}}}(\alpha_k + (\ell - 1 - n)\alpha_{i_{j+1}})$,
\cite[Theorem 9.3]{decpro3}. Therefore 
\[
T_{w_j}(E_kE_{i_{j+1}}^{\ell - 1})\otimes 1 = \sum_{n=0}^{\ell -
1}T_{w_j}(E_{i_{j+1}}^{n})T_{w_{j+1}}(u_n') \otimes 1 = 0,
\]
as required.
\end{proof}
\subsection{Proof of Theorem \ref{simpqg}} Suppose that $V_{\chi}(t)$ is simple, where $\chi = \chi_u\chi_s$ lies over $g = g_ug_s \in G$. We retain the notation of (\ref{simpqg}), but in view of (\ref{quanred}) we can and do assume that $g$ is exceptional. Applying Lemma \ref{exceptional} and using the notation introduced there, as $g=g_ug_s$ is exceptional we can assume without loss of generality
that $g_s = s_m$ for some $m$, $0\leq m\leq r$. Hence, by Lemma
\ref{exceptional}.2, a set of simple roots for $C_G(g_s)$ is given by $\Delta$ if
$m=0$ or by $\{-\alpha_1,\ldots , -\widehat{\alpha_m},\ldots ,-\alpha_r ,
-\beta_m\}$ if $1\leq m\leq r$.

Suppose the theorem is false. Then we can pick a component $\Phi_i$
where $g_i$ is irregular and $\alpha(t)^2 \neq \ep^{-2(\rho ,\alpha)}$
for some $\alpha \in \Delta_i$. By \cite[Proposition 4.1 and Theorem
4.9]{humCC} we can assume, without loss of generality, that there is
an element $\beta\in \Delta_i$ such that the $U_{\beta}$-component of
$g_u$ is trivial. Let $\beta_1 ,\ldots , \beta_s$ be a minimal length connected chain in
$\Delta_i$ such that the $U_{\beta_1}$-component of $g_i$ is trivial
and $\beta_s(t)^2\neq \ep^{-2(\rho, \beta_s)}$. 

Using the notation of Section (\ref{quanbab}), consider the element
$s_{\beta_s}{\scriptscriptstyle \bullet} t \in T$. By construction
\begin{equation}
\label{eq:rambofirstblood}
\beta_{s-1}(s_{\beta_s}{\scriptscriptstyle \bullet} t)^2 = \ep^{2(\rho,
-<\beta_{s-1},\check{\beta_s}>\beta_s)}\beta_{s-1}(t)^2\beta_s(t)^
{-2<\beta_{s-1},\check{\beta_s}>}.
\end{equation}
By minimality $\beta_{s-1}(t)^2 = \ep^{-2(\rho, \beta_{s-1})}$. It
is easy to check that $<\beta_{s-1},\check{\beta_s}>$ is
invertible modulo $\ell$, so it follows from (\ref{eq:rambofirstblood}) that
$\beta_{s-1}(s_{\beta_s}{\scriptscriptstyle \bullet} t)^2 \neq \ep^{-2(\rho,
\beta_{s-1})}$. By (\ref{eq:quantlinkage}) and Theorem \ref{quanbab} the modules $V_{\chi}(t)$ and
$V_{\chi}(s_{\beta_s}{\scriptscriptstyle \bullet} t)$ are in the same block, but as
$V_{\chi}(t)$ is simple and has dimension equal to the PI-degree of $\UE$ its anihilator in $Z(\UE)$ must belong to the Azumaya locus. Hence $V_{\chi}(t)$ and $V_{\chi}(s_{\beta_s}{\scriptscriptstyle \bullet} t)$ must, in fact,
be isomorphic. Thus we can assume, without loss of generality, that
there exists $\beta\in \Delta_i$ such that the $U_{\beta}$-component
of $g_i$ is trivial and $\beta(t)^2\neq \ep^{-2(\rho,\beta)}$.  

Define $w = s_{i_1}\ldots s_{i_v} \in W$ and $\alpha\in \Delta$ (the simple roots of $\g$) as
follows: if $-\beta\in \Delta$ let $w=e$ and $\alpha = -\beta$; if
$\beta=-\beta_m$ then let $w = w_m$ and $\alpha$ be determined as in Lemma \ref{calculation}. 
Let $w'=s_{i_1}\ldots s_{i_v}s_{\alpha}$ and for $1\leq j\leq
v+1$ let $\gamma_j = s_{i_1}\ldots s_{i_{j-1}}\alpha_{i_j}$ (so
$\gamma_1 = \alpha_{i_1}$ and $\gamma_{v+1}=\beta$).
We claim that $\chi_u(F_{\gamma_j}^{\ell}) = 0$ for all $j$, $1\leq
j\leq v+1$. Indeed if $j \leq v$ then by Lemmas \ref{calculation}.2 and \ref{calculation}.3  $\gamma_j$ is not a
root of $C_G(g_s)$ which means that $\chi_s(K_{\gamma_j}^{2\ell})=
\gamma_j(s_m)\neq 1$. By Lemma \ref{hearmysong} this implies that
$\chi_u(F_{\gamma_j}^{\ell}) = 0$. For the case $\gamma_t$ the
hypothesis that $g_u$ has a trivial $U_{-\gamma_{v+1}}$ component implies,
as in the proof of Lemma \ref{hearmysong}, that $\chi_u(F_{\gamma_{v+1}}^{\ell})=0$.

We can iterate Lemma \ref{twist}, using the fact that $V_{\chi}(t)$ is simple
to obtain an isomorphism
\[
V_{\chi}(t) \cong V_{\chi}^w(t'),
\]
where $t'(K_{\lambda}) = \ep^{(\lambda, \gamma_1 +\ldots
+\gamma_{t})}t(K_{\lambda})$. Let $\tilde{R}$ be the algebra of
$\UE$ generated by $U_{\ep}^{\geq 0}$ and the element $F_{\alpha}$ and
let $R$ be the quotient obtained factoring out the ideal generated by
$F_{\alpha}^{\ell}$, the elements $E_{\beta}^{\ell}$ for all $\beta\in
\Phi^+$ and the elements
$K_{\lambda}^{\ell}-(w^{-1}\chi_s)(K_{\lambda}^{\ell})$ for all
$\lambda \in X$. Then there is an
isomorphism
\[
V_{\chi}^w(t') \cong U_{\chi}(\g)
\otimes_{T_w(R)} T_w(R) \otimes_{T_w(U_{w^{-1}\chi_s}^{\geq 0})}
k_{t'}.
\]
The isomorphism from $R\otimes_{U_{w^{-1}\chi_s}^{\geq
0}}k_{w^{-1}(t')}$ to
$T_w(R)\otimes_{T_w(U_{w^{-1}\chi_s}^{\geq 0})} k_{t'}$ given by sending $r\otimes 1$ to $T_w(r)\otimes 1$
transforms $R$-modules into $T_w(R)$-modules. 

We claim that $R\otimes_{U_{w^{-1}\chi_s}^{\geq
0}}k_{w^{-1}(t')}$ has a non-trivial
quotient. Indeed $\tilde{R}$ is just a parabolic extension of
$U_{\ep}(\mathfrak{sl_2})$ so, by \cite[Proposition 2.12]{janqg}, $R\otimes_{U_{w^{-1}\chi_s}^{\geq
0}}k_{w^{-1}(t')}$ is simple if and only if $\alpha(w^{-1}(t'))^2 =
\ep^{2(\rho ,\alpha)}$, since $w^{-1}\chi_s(K_{\alpha}^{2\ell}) =
\chi_s(K_{\beta}^{2\ell}) = 1$. The equality $\rho - w\rho = \gamma_1
+\ldots +\gamma_{v}$, however, shows that
\[
\alpha(w^{-1}(t'))^2 = \beta(t')^2 =
 \ep^{2(w^{-1}(\rho-w\rho),\alpha)}\beta(t)^2.
\]
Hence $\alpha(w^{-1}(t'))^2=\ep^{-2(\rho,\alpha)}$ if and only if
$\beta(t)^2 = \ep^{-2(\rho , \beta)}$. 

By hypothesis
$\beta(t)^2\neq \ep^{-2(\rho,\beta)}$, so $R\otimes_{U_{w^{-1}\chi}^{\geq
0}}k_{w^{-1}(t')}$ has a non-trivial quotient
and therefore so too does $V_{\chi}(t)$, contradicting the assumed simplicity.
\subsection{}
We will now describe the fully Azumaya locus. Recall that we say
$\chi\in \mathrm{Maxspec}(Z_0)$ is regular if it lies over a regular element of $G$.
\begin{thm}
Let $\chi\in \mathrm{Maxspec} (Z_0)$. Then $\mathfrak{m}_{\chi}$ is in the fully Azumaya locus - that is, the simple modules of $U_{\chi}(\g)$ all have
dimension $\ell^N$ - if and only if $\chi$ is regular.
\end{thm}
\begin{proof}
Let $\chi$ lie over $g=g_sg_u$, as usual. Note that $\Phi'=\emptyset$
if and only if $C_G(g_s)=T$ if and only if $g=g_s$ is regular
semisimple. 

Suppose all simple $U_{\chi}(\g)$-modules have
dimension $\ell^N$. Then by Theorem \ref{simpqg} either $\pi(\chi_u)$ is
regular in $C_{G'}(g_s)=C_G(g_s)$ or $\Phi'=\emptyset$. In both cases
this implies that $g$ is regular, as required.
\end{proof}
\begin{rem}
Suppose that $\ell$ is prime to the index of connection of $\g$. Then
the above theorem remains true for any $Q\subset M\subset X$. For
under this restriction simple baby Vermas over $X$ restrict to simple
baby Vermas over $M$. Indeed let $0\neq v \in V_{\chi}(t)$, a simple baby
Verma over $X$. We will show that $U_{\chi, M}(\g)v = V_{\chi}(t)$. As
$U_{\chi}^+$ acts nilpotently we can assume without loss of generality
that $U_{\chi}^+v=0$. After diagonalisation we may also assume that
$v$ is an eigenvector for the action of $K_{\lambda}$ with $\lambda\in
M$. By hypothesis, the index of connection is invertible hence the
diagonalisation for $K_{\mu}$ with $\mu\in X$ is identical. Since
$V_{\chi}(t)$ is a simple $U_{\chi}(\g)$-module we therefore see that
$U^-v = V_{\chi}(t)$ as required. 
\end{rem}
\subsection{}
The analysis of section (\ref{reduce}) also allows one completely to describe the
unramified maximal ideals of $Z$. We will use the notation of sections
(\ref{exceptional}) and (\ref{reduce}). So assume that $\chi=\chi_u\chi_s\in \mathrm{Maxspec}( Z_0)$ lies
over $g=g_ug_s$. In particular recall from (\ref{eq:barbaracartland}) the definition of
$\Phi'$ and set
$
X' = \mathbb{Z}\{ \varpi_i : \alpha_i\in \Phi'\}$.
Let $\tilde{p}:T\longrightarrow T'$ be the morphism associated with
the inclusion $\C [X']\subseteq \C [X]$. Let $s=\dim X-\dim
X'$. Finally let $\tilde{\Phi}'\subseteq \Phi'$ be the set $\{
\beta\in \Phi : \beta(g_s)=1\}$. 
\begin{prop}
Keep the above notation. The unramified maximal
ideals of $Z$ lying over $\mathfrak{m}_{\chi}$ are of the form
$\mathrm{Ann}_Z(V_{\chi}(t))$ where
$\alpha((\tilde{p})(t))^2=\ep^{-2(\rho,\alpha)}$ for all $\alpha \in
\tilde{\Phi}'$. In particular if $\ell$ is prime to the index of connection of
$\Phi'$ then there are exactly $\ell^s$ such ideals.
\end{prop}
\begin{proof}
The first statement is an immediate consequence of Theorem \ref{projqg} and
(\ref{eq:caddyshack}). For the second statement note first that if we know the
value of $\alpha(t')$ for all $\alpha \in \tilde{\Phi}'$ then we know
the value of $\alpha(t')$ for all $\alpha \in \Phi'$ (since, by section
(\ref{hearmysong}),
$\tilde{\Phi}'$ has, up to conjugacy, a basis of maximal rank
consisting of simple roots of $\Phi'$ and its longest root). If $\ell$
is prime to the index of connection of $\Phi'$ then we can calculate
the values of $\varpi (t')$ given the values of $\alpha_i(t')$ for
$\alpha_i\in \Phi'$ which means that $t'$ is completely determined in
this case. Hence there is only one possible value for $\tilde{p}(t)$
so there are only $\ell^s$ preimages. Suppose that
$\mathrm{Ann}_Z(V_{\chi}(t))=\mathrm{Ann}_Z(V_{\chi}(u))$ for two such preimages. Then,
by (\ref{eq:quantlinkage}), $u=w{\scriptscriptstyle \bullet} t$ for some $w\in W$. As $\chi_s = u^{\ell} =
(w{\scriptscriptstyle \bullet} t)^{\ell} = w(t^{\ell}) = w\chi_s$ we see that $w$ belongs
to the Weyl group of $C_G(g_s)$. Since both $t$ and $u$ lie over
$t'\in T'$ as above it follows that they are fixed under the dot
action of the Weyl group of $C_G(g_s)$, implying that $t=u$ as required.
\end{proof}
\newpage
\appendix
\section{}
Here we list the elements of the Weyl group and the simple roots
required for Lemma \ref{calculation}.

\begin{tabular}{|l|l|l|l|} \hline
Type & $m$ & $w_m$ & $\alpha^m$ \\ \hline
$A_r$ & $1\leq m\leq r$ & $e$ & $\alpha_m$ \\ \hline
$B_r$ & 1 & $e$ & $\alpha_1$ \\
& $1< m \leq r$ & $s_ms_{m+1}s_{m+2}\cdots s_{r-1}s_rs_{r-1}\cdots
s_{m+1}s_m$ & $\alpha_{m-1}$ \\ \hline
$C_r$ & $1\leq m\leq r$ & $ s_ms_{m+1}\cdots s_{r-1} $ & $\alpha_r$ \\
\hline
$D_r$ & $1,r-1$ or $r$ & $e$ & $\alpha_1 , \alpha_{r-1} , \alpha_r
$ \\
 & $2\leq m \leq r-2$ & $s_ms_{m+1}\cdots
s_{r-2}s_rs_{r-1}s_{r-2}\cdots s_{m+2}s_{m+1}s_{m-1}$ & $\alpha_m$ \\
\hline
$F_4$ & 1 & $s_1s_2s_3s_2s_4s_3s_2$ & $\alpha_1$ \\ 
& 2 & $s_2s_3s_2s_1s_4s_3$ & $\alpha_2$ \\
 & 3 & $s_3s_2s_1s_4s_3$ & $\alpha_2$ \\
& 4& $s_4s_3$ & $\alpha_2$ \\ \hline
$G_2$ & 1 & $s_1$ & $\alpha_2$ \\
& 2 & $s_2s_1 $ & $\alpha_2$ \\ \hline   
$E_6$ &1& $e$ & $\alpha_1$ \\
&2& $s_2s_4s_5s_6s_3s_1s_4s_3s_5s_4$ & $\alpha_2$ \\
&3& $s_3s_1s_4s_5s_2s_4$ & $\alpha_3$ \\
&4& $s_4s_5s_6s_3s_1s_4s_3s_5s_2$ & $\alpha_4$ \\
&5& $s_5s_6s_4s_3s_2s_4$ & $\alpha_2$ \\
&6& $e$ & $\alpha_6$ \\ \hline
$E_7$ & 1 &  $s_1s_3s_4s_2s_5s_4s_3s_6s_5s_4s_1s_2s_3s_4s_5s_6$ &
$\alpha_7$ \\ 
&2& $s_2s_4s_5s_6s_3s_1s_4s_3s_5s_4$ & $\alpha_2$ \\
&3& $s_3s_4s_2s_5s_4s_3s_6s_5s_4s_1s_2s_3s_4s_5s_6$ & $\alpha_7$ \\
&4& $s_4s_2s_5s_4s_3s_6s_5s_4s_1s_2s_3s_4s_5s_6$ & $\alpha_7$ \\
&5& $s_5s_4s_3s_6s_5s_4s_1s_2s_3s_4s_5s_6$ & $\alpha_7$ \\
&6& $s_6s_5s_4s_2s_3s_4s_5s_6$ & $\alpha_7$ \\
&7& $e$ & $\alpha_7$ \\ \hline
$E_8$ &1& $s_1s_3s_4s_2s_5s_4s_3s_6s_5s_4s_1s_2s_3s_4s_5s_6$ &
$\alpha_7$ \\ 
&2&
$s_2s_4s_3s_5s_4s_2s_6s_5s_4s_3s_7s_6s_5s_4s_1s_2s_3s_4s_5s_6s_7$ &
$\alpha_8$ \\
&3&
$s_3s_1s_4s_3s_5s_4s_2s_6s_5s_4s_3s_7s_6s_5s_4s_1s_2s_3s_4s_5s_6s_7$ &
$\alpha_8$ \\
 &4&
$s_4s_2s_3s_1s_4s_3s_5s_4s_2s_6s_5s_4s_3s_7s_6s_5s_4s_1s_2s_3s_4s_5s_6s_7$
& $\alpha_8$ \\
& 5 &
$s_5s_4s_2s_3s_1s_4s_3s_5s_4s_2s_6s_5s_4s_3s_7s_6s_5s_4s_1s_2s_3s_4s_5s_6s_7$
& $\alpha_8 $ \\
& 6 &
$s_6s_5s_4s_2s_3s_1s_4s_3s_5s_4s_2s_6s_5s_4s_3s_7s_6s_5s_4s_1s_2s_3s_4s_5s_6s_7$
& $\alpha_8$ \\
& 7 &
$s_7s_6s_5s_4s_2s_3s_1s_4s_3s_5s_4s_2s_6s_5s_4s_3s_7s_6s_5s_4s_1s_2s_3s_4s_5s_6s_7$
& $\alpha_8$ \\
& 8 &
$s_8s_7s_6s_5s_4s_2s_3s_1s_4s_3s_5s_4s_2s_6s_5s_4s_3s_7s_6s_5s_4s_1s_2s_3s_4s_5s_6s_7$
& $\alpha_8$ \\ \hline
\end{tabular}

\end{document}